%% file: main.tex
\documentclass[reqno]{amsart}
\usepackage[margin=1.4in]{geometry}
\usepackage[]{amssymb} 
\usepackage[]{amsthm} 
\usepackage{amsmath}
\usepackage{mathrsfs}
\usepackage{mathtools}
\usepackage[hidelinks]{hyperref}
\usepackage{enumerate}
\usepackage{array}
\usepackage{float}
\usepackage{arydshln} 
\usepackage{bbm}
\usepackage{cleveref}
\usepackage{scalerel,stackengine}
\usepackage[dvipsnames]{xcolor} 
\usepackage{bm}

\usepackage{autonum}
\newtheorem{theorem}{Theorem}[section]
\newtheorem*{theorem*}{Theorem}
\newtheorem{prop}{Proposition}[section]
\newtheorem*{prop*}{Proposition}
\newtheorem{conj}{Conjecture}
\newtheorem{lemma}{Lemma}[section]

\theoremstyle{definition}

\theoremstyle{remark}
\newtheorem{remark}[theorem]{Remark}
\newtheorem*{remark*}{Remark}

\numberwithin{equation}{section}

\newcommand{\brackets}[1]{\left(#1\right)}  % brackets

\newcommand{\Z}{\mathbb{Z}}
\renewcommand{\P}{\mathbb{P}}
\newcommand{\N}{\mathbb{N}}

\renewcommand{\phi}{\varphi}

\newcommand{\PP}{\mathbb{P}}

\newcommand{\mA}{\mathcal{A}}
\newcommand{\mB}{\mathcal{B}}

\newcommand{\mFb}{\mathcal{F}_{b,\eta}}
\begin{document}

\title[Shifted Prime Divisors]{Integers divisible by a shifted prime in a given interval}

%    Information for authors
\author[Abi Abdallah]{Rebecca Abi Abdallah} 
\author[Kovaleva]{Valeriya Kovaleva}
\author[Schlitt]{Jeremy Schlitt}
\author[Tardy]{N\'eo Tardy}
\address{D\'epartement de Math\'ematiques et de Statistique, Universit\'e de Montr\'eal, CP 6128 succ. centre, Centre Ville, Montr\'eal, Qc 3HJ3J7, CANADA}
\thanks{}

%    General info
\subjclass[2020]{}

\date{\today}

\begin{abstract}
        In this paper we study the behaviour of $H^*(x,y,z)$, the number of integers less than $x$ possessing a divisor in the interval $(y,z]$ of the form $p-1$, where $p$ is a prime, for all values of $y = y(x)$ and $z= z(y)$. We observe multiple phase transitions at critical values of $z$ in terms of $x$ and $y$ guided largely by the anatomy of $n$. Our results generalize a result of Ford from 2017, which corresponds to the case $z = x$. 
\end{abstract}

\maketitle

\section{Introduction}
\label{sec:intro}
\input{01_new_intro}

\subsection{Outline of the proof}
\label{sec:structure}

\input{02_technical_details}

\section{Preliminaries}
\label{sec:prelim}

\input{03_preliminaries}
% \newpage

\section{Asymptotic for short intervals}

\label{sec:short}
\input{04_short_intervals}

\section{Order of magnitude for longer intervals} 
\label{sec:long}
\input{05_long_intervals}

\bibliographystyle{plain}
\bibliography{main.bib}{}
\end{document}

%% file: 01_new_intro.tex
Let $\tau(n;y,z) = \sum_{d \mid n,\ y < d \le z} 1$ be the number of divisors of $n$ in the interval $(y,z]$. Understanding the statistics of $\tau(n;y,z)$ is one of the central problems in the anatomy of integers and the statistics of product sets. 

Denote the number of integers that have at least one divisor in $(y,z]$ by
\begin{equation}
    H(x,y,z) = \#\{n \le x: \tau(n;y,z) \ge 1\}\,.
\end{equation}
The quantity $H(x,y,z)$ with $z = 2y$ was first studied by Erd\H{o}s \cite{erdHos1955some,erdos1960ob} in order to count distinct entries in the $N \times N$ multiplication table $\mathcal{A}(N)$. In particular, Erd\H{o}s obtained the first non-trivial bound $|\mathcal{A}(N)| = o(N^2)$ by showing that $H(x,\sqrt{x}/2,\sqrt{x}) = o(x)$. Inspired by Erd\H{o}s and newly found applications (see e.g. \cite{cobeli2003jumping,erdos1981structure,tenenbaum1976deux,heath2003linear,ford2005maximal}), several authors worked towards determining the order of $H(x,y,z)$ for all $(y,z]$. Turns out, the behaviour of $H(x,y,z)$ depends significantly on the anatomy of integers producing the divisor and as a result it exhibits multiple phase transitions as a function of $\log(z/y)$, the logarithmic length of the interval. We particularly note the efforts of Hall and Tenenbaum \cite{Hall_Tenenbaum_1988}, who established almost tight though not exact upper and lower bounds on $H(x,y,z)$ in many cases. 

For $Cy \le  z \le x^{0.99}$, it was ultimately Ford \cite{ford2008distribution} who established the exact order of magnitude of $H(x,y,z)$ in his breakthrough work using a ballot theorem-type argument. Most notably, Ford showed that 
\begin{equation}\label{thm:FORDANNALS0}
    H (x,y,z)\asymp xu^{\delta}\log (2/u)^{-3/2}\,,
\end{equation}
where $u = \frac{\log (z/y)}{\log z}$, and $\delta = 1 - \frac{1+\log\log 2}{\log 2}$ is the Erd\H{o}s–Ford–Tenenbaum constant, or the multiplication table constant. That is, when $z = 2y \le x^{0.99}$
\begin{equation}\label{thm:FORDANNALS}
    H (x,y,2y) \asymp x (\log y)^{-\delta}(\log\log y)^{-3/2}
\end{equation}
is smaller than the heuristic estimate by a factor of $\log \log y$. This unusual extra factor is coming from the clustering of divisors of $n$ on the logarithmic scale caused by irregularities in the distribution of prime factors of $n$, a profound structural phenomenon observed by Ford.

Whether an asymptotic $H(x,y,z) = x \rho(y,z) (1+ o_{x \to \infty}(1))$
exists for all $y,z$ remains a difficult open question, however, in their recent work, Green and Sawhney \cite{green2026proportion} showed that the permutation analog of $H(x,y,2y)$ indeed possesses an asymptotic when $y = x^{o(1)}$.

To expand upon the structural study of divisors, one can restrict the integers in question to a special set $\mA$. For example, one can count generic divisors $d$ of integers $n \in F(\Z)$ for $F$ an irreducible polynomial over $\Z$ \cite{tenenbaum1990question}, or $n$ in an arithemtic progression \cite{ford2005maximal}, or when $n$ is a rough integer \cite{ford2021rough}. More recently, the third author \cite{schlitt2026multiplication} extended the latter by considering integers having all their prime factors in a given set of prescribed density. Another interesting example is when $n = p - 1$ is a shifted prime, see \cite{pappalardi1996order,erdos1999order,indlekofer2002divisors, ford2008distribution, koukoulopoulos2010generalized} as well as \cite{ koukoulopoulos2010generalized, koukoulopoulos2010divisors}.

On the other hand, one could also count divisors $d\in\mA$ of generic integers $n$. In both cases, it is interesting to ask whether the clustering phenomenon mentioned above persists, and if one could obtain more precise information about the asymptotic behaviour of $H(x,y,z)$. The case $d\in\mA$ appears to have received less attention, and is the subject of the present paper. 

\subsection{Shifted prime divisors}
Denote the number of integers up to $x$ divisible by at least one shifted prime $p-1$ such that $p \in (y,z]$ by
\begin{equation}
    H^*(x,y,z) := \#\{n \le x: \exists p \in (y,z],\ p-1 \mid n\}\,.
\end{equation}
The first instance of such a function (with $z=x$) is due to Erd\H{o}s and Wagstaff \cite{erdHos1980fractional}, who showed that $H^*(x,y,x) \ll x (\log y)^{-c}$ for some constant $c>0$. Though studying $H^*(x,y,x)$ is thematically a question about the distribution of divisors of an integer, it was initially motivated by a completely different application, namely, the study of fractional parts of the Bernoulli numbers. Subsequently, McNew, Pollack and Pomerance \cite{mcnew2017numbers} determined that one can take $c=\delta +o(1)$, and proved a matching lower bound. Ford \cite{ford2017integers} resolved this question for all intervals $[y,x)$. In particular, he showed that for $3 \le y\leq x^{0.99}$ 
\begin{equation}\label{thm:FORDSHIFTED}
            H^*(x,y,x) \asymp x (\log y)^{-\delta} (\log \log y)^{-1/2}\,.
        \end{equation}
Ford's proof also implies the same order-of-magnitude estimate for $H^*(x,y,z)$ when $z \geq y^{1+\varepsilon}$ for small fixed $\varepsilon>0$. 
Contrary to \eqref{thm:FORDANNALS0}, the order of magnitude in \eqref{thm:FORDSHIFTED} does not exhibit any unusual factors. In other words, $H^*(x,y,x)$ is not affected by the clustering of divisors. 

In this paper we establish the order of magnitude of $H^*(x,y,z)$ for almost all intervals $[y,z)$ generalizing the work of Ford. Curiously, this result can simultaneously be seen as an extension of two somewhat independent problems \eqref{thm:FORDANNALS0} and \eqref{thm:FORDSHIFTED}.

\begin{theorem}\label{mainthm} Let $3 \leq y < z \leq x$, and let $z^\prime = \min\{z,x/y\} \to \infty$ as $x \to \infty$. Let $\varepsilon,\eta \in (0,1/10)$ be small and fixed. Let $ Q(\rho) = \rho \log (\rho/e) + 1$ and $\gamma_0 = Q(2) = \log 4 -1$. Set
\[
u = \frac{\log(z/y)}{\log z},\quad u^\prime = \frac{\log(z/y)}{\log z^\prime}, \quad \gamma = \frac{\log (1/u)}{\log \log z^\prime}\,
\]
such that $y = z^{1-u}$ and $u = (\log z^\prime)^{-\gamma}$. Then the following holds.
\begin{enumerate}[(i)]

\item If $u \le (\log z^\prime)^{-(\gamma_0+\eta)}$, or equivalently $\gamma > \gamma_0 + \eta$, then for $y \ge z - z^{17/30+\varepsilon}$ 
\begin{equation}
   x^{-1} H^*(x,y,z) = u (1+ o_{\varepsilon,\eta}(1))\,.
\end{equation}

\item If $\gamma \in (\eta,\gamma_0 - \eta)$, then
\begin{equation}
        x^{-1} H^*(x,y,z) \asymp_\eta (\log z^\prime)^{-Q(\frac{1+\gamma}{\log 2})} (\log \log z^\prime)^{-1/2}\,.
        \end{equation}

 \item (First transition) If $|\gamma-\gamma_0| \le 0.01$, then
\begin{equation}
    \begin{aligned}
        x^{-1}H^*(x,y,z) \asymp \begin{cases}
            u, & \text{if $\gamma \ge \gamma_0$}\,,\\
                (\log z^\prime)^{ -Q(\frac{1+\gamma}{\log 2})} ((\gamma_0-\gamma)\sqrt{\log\log z^\prime} + 1)^{-1}, & \text{otherwise.}
        \end{cases}
    \end{aligned}
\end{equation}
\item (Second transition) Let $\lambda_0 = \log (1/u^\prime)$. If $0 \le \gamma \le 0.01$, then 
                \begin{equation}
          (\log z^\prime)^{-Q(\frac{1+\gamma}{\log 2})} (\log \log z^\prime)^{-1/2} V^{-1}  \ll x^{-1}H^*(x,y,z) \ll (\log z^\prime)^{-Q(\frac{1+\gamma}{\log 2})} (\log \log z^\prime)^{-1/2}\,,
\end{equation}
where $V = V(\gamma,\lambda_0) \le \exp\left(C(\varepsilon^\prime) \min\{1/\gamma^2, 1+\lambda_0\}^{1/2+\varepsilon^\prime}\right)$ for any $\varepsilon^\prime > 0$.
\end{enumerate}
\end{theorem}
\begin{remark}
   In \Cref{mainthm}(i) if $y \le z - z^{1-\varepsilon}$ it may occur that there are no primes in the interval $(y,z]$ so that $H^*(x,y,z) = 0$. For such short intervals our argument does not rely on the knowing the exact count and instead shows that 
            \begin{equation}
            x^{-1}H^*(x,y,z) = \left( \sum_{y< p \le z}{ \frac{1}{p-1}}\right)(1+ o(1))\,.
        \end{equation}
In this case an asymptotic for $H^*(x,y,z)$ is achievable in the range depending on existing estimates for primes in short intervals; the recent result of Guth and Maynard \cite{guth2026new} allows us to take $y$ as small as $ y \ge z- z^{17/30+\varepsilon}$. 
\end{remark}
\begin{remark}
    Ford's result \cite{ford2017integers} in full generality corresponds to $x=z$. This case is included in our Theorem as $\gamma = \frac{1}{\alpha} - 1$ ($\alpha = \frac{\log\log z^\prime}{\log\log z}$ is the same as in the notation of \cite{ford2017integers}) and $\lambda_0 = 0$. In particular,  \Cref{thm:FORDSHIFTED} corresponds to $\gamma = 0$, in which case $Q(\frac{1+\gamma}{\log 2}) = Q(\frac{1}{\log 2}) = \delta$.
\end{remark}

Our Theorem produces the exact order of magnitude in all cases but one: when $V \to \infty$ in \Cref{mainthm}(iv). If $\gamma \gg 1$ or $\lambda_0 \ll 1$, then $V \ll 1$ poses no issue. When $\gamma = \gamma(u,z^\prime) \to 0$ and $\lambda_0 = \lambda_0(z/y,z^\prime) \to \infty$ simultaneously (this does not include $\gamma \ll (\log\log z^\prime)^{-1}$ or $x=z$) we have $\min\{1/\gamma^2, \lambda_0 + 1\} \le \lambda_0$ and thus
\begin{equation}
    V^{-1} \ge \brackets{\frac{\log z/y}{\log z^\prime}}^{o(1)} \ge (\log z^\prime)^{-o(1)}
\end{equation}
limits our result to 
\begin{equation}
    H^*(x,y,z) = x (\log z^\prime)^{-Q(\frac{1+\gamma}{\log 2}) + o(1)}\,.
\end{equation}
One may ask whether this failure to produce the correct order is the failure of the method, or the effect of clustering as in \eqref{thm:FORDANNALS0}. We are inclined to believe the latter meaning that neither our upper nor our lower bound are sharp in this case. We conjecture that the true order of magnitude is as follows.

\begin{conj} Let $\lambda_0 = \log\log z^\prime-\log\log (z/y)$. Then for $\gamma \le 0.01$ we have
\[
x^{-1}H^*(x,y,z) \asymp (\log z^\prime)^{-Q(\frac{1+\gamma}{\log 2})} (\log\log z^\prime)^{-1/2} \cdot (\gamma + (\lambda_0+1)^{-1/2})\,.
\]
\end{conj}

Let us briefly discuss the main features of our result and explain the origin of the first transition at $\gamma = \gamma_0$ and the conjectured second transition at $\gamma = 0$ in broad terms. Technical details as well as the outline of the proof can be found in \S\ref{sec:structure}.

\subsection{Discussion of the main theorem}
If $z^\prime = \min\{z,x/y\}$, on average $\tau(n;y,z)$ depends on $\Omega(n,z^\prime)$, the number of prime factors of $n$ less than $z^\prime$, and is proportional to the logarithmic length of the interval. That is, an integer with $\Omega(n,z^\prime)=k$ has on average
\begin{equation}\label{eq:divisors_heuristic0}
   \frac{\#\{n = dm \le x: y < d \le z,\, \Omega(n,z^\prime) = k \}}{\#\{n \le x: \Omega(n,z^\prime) = k \}} \asymp 2^{k} \cdot \frac{\log (z/y) }{\log z^\prime}
\end{equation}
divisors in the interval $(y,z]$. Quantity \eqref{eq:divisors_heuristic0}, while providing a valid upper bound for the proportion of integers with a divisor in an interval, is not necessarily a good estimate: at the very least, it needs to beat the trivial bound of $1$. When $z/y \to \infty$, it is further inflated by $\tau(n_{\le z/y})$, the number of $(z/y)$-smooth divisors of $n$. 

The factor $2^k$ in \eqref{eq:divisors_heuristic0} corresponds to the number of ways to distribute $k$ prime factors of $n = dm$ between two sets. Another way one can interpret $2^k$ is as follows. We know that $\Omega(n,z^\prime)$ behaves roughly like a Poisson random variable with parameter $\lambda = \log\log z^\prime$ \cite{Hall_Tenenbaum_1988}. If $n = dm$ is a product of two integers, we have that $P_1 = \Omega(d,z^\prime) \approx P_2 = \Omega(m,z^\prime) \approx k/2$ are also approximately independent Poisson random variables with parameter $\lambda = \log\log z^\prime$. Then the representation function $\tau(n;y,z)$ counting pairs $(d,m)$ ends up being the distribution function of $P_1+P_2$, which is itself a Poisson random variable with parameter $2\lambda$.

Similarly for shifted primes, an integer $n$ with $k$ prime factors up to $z^\prime$ has on average  
\begin{equation}\label{eq:shifted_prime_heuristic0}
    \frac{\#\{n = (p-1)m \le x: y < p \le z,\, \Omega(n,z^\prime) = k \}}{\#\{n \le x: \Omega(n,z^\prime) = k \}} \asymp  2^{k} \cdot \frac{\log z/y}{\log z^\prime} \cdot \frac{1}{\log z} 
    = \frac{2^{k}u}{\log z^\prime}
\end{equation}
shifted prime divisors in the interval $(y,z]$ producing a function likewise proportional to the logarithmic length of the interval. In this case, one expects that $n$ needs to have at least $\log z$ divisors in the interval $(y,z]$ in order to have at least one shifted prime divisor. 

In both cases the interaction between the two Poisson random variables with parameters $\log\log z^\prime$ and $2\log\log z^\prime$ is qualitatively different depending on the logarithmic length of $(y,z]$. Ignoring the issues with the $(z/y)$-smooth part for now, we observe that \eqref{eq:shifted_prime_heuristic0} beats the trivial bound of $1$ when
\[
\Omega (n,z^\prime) \le \frac{1}{\log 2} \log (u^{-1} \log z^\prime) = 
\frac{1+\gamma}{\log 2} \log\log z^\prime\,.
\]
At the same time, the Poisson distribution with parameter $\lambda^\prime = 2 \log \log z^\prime$ concentrates on
\begin{equation}\label{eq:Omeganz_i 0}
    \Omega (n,z^\prime) \le 2 \log\log z^\prime + O((\log\log z^\prime)^{1/2+\varepsilon})\,.
\end{equation}
Thus, one needs to consider $n$ with
\begin{equation}
    \Omega (n,z^\prime) \lesssim \min\left\{ \frac{1+\gamma}{\log 2},2 \right\}\log\log z^\prime\,.
\end{equation}
This observation marks a transition in the anatomy of integers with the largest contribution to $H^*(x,y,z)$ at $\frac{1+\gamma}{\log 2} = 2$, or $\gamma = \log 4 - 1 = \gamma_0$. When $\gamma > \gamma_0$, estimate \eqref{eq:shifted_prime_heuristic0} is almost always accurate meaning that the interval $(y,z]$ most of the time has a single shifted prime divisor if it has one at all, which allows up to establish an asymptotic for $H^*(x,y,z)$. For $\gamma \lesssim \gamma_0$ this is not the case anymore, so we are only able to establish the order of magnitude. 

The second transition occurs near $\gamma=0$. In this regime the main contribution to $H^*(x,y,z)$ comes from $n$ with
\begin{equation}
     \Omega (n,z^\prime) = \frac{1}{\log 2} \log\log z^\prime + o(\log\log z^\prime)\,,
\end{equation}
which is close to the ``critical'' number of prime factors in \cite{ford2008distribution} causing $\{\log d\}_{d \mid n}$ to cluster. Similarly to \cite{ford2008distribution}, we expect that our bounds in \Cref{mainthm}(iv) are not sharp due to clustering of $\{\log (p-1)\}_{p-1\mid n}$ caused by irregularities in the distribution of prime factors of $n$. However, for this phenomenon to affect $H^*(x,y,z)$ in contrast with \eqref{thm:FORDSHIFTED}, the interval $(y,z]$ has to be logarithmically short, that is, it requires
\begin{equation}
    u^\prime = \frac{\log (z/y)}{\log z^\prime} = e^{-\lambda_0} = o (1)\,, 
\end{equation}
and thus $\lambda_0 \to \infty$. The only regime when the two conditions $\gamma \to 0$ and $\lambda_0 \to \infty$ co-occur happens to be relatively brief.

\subsection*{Acknowledgements}
We thank Dimitris Koukoulopoulos for his guidance throughout the early stages of this project. We also thank Cihan Sabuncu for giving a helpful lecture on sieve methods in summer 2024.

\noindent 
The project was supported by NSERC (grant No. 2018-05699). VK and JS were additionally supported by Courtois Chaire en recherche fondamentale II. RAA and JS were supported by NSERC (grant No. 2024-05850). RAA was supported by the CRM-ISM undergraduate research scholarship and FRQ-PBEEE (grant No. 370608). VK was supported by FRQNT (grant No. 373840 and 345672). JS was supported in part by the University of Montr\'eal's Bourse d'Excellence - Banque Nationale. NT was supported by ENS Paris-Saclay.

\subsection*{Notation} 
We adopt the standard Vinogradov notation $\ll,\gg,\sim,\asymp$ as well as Big-Oh notation throughout the paper. The function $\log x$ stands for the natural logarithm of $x$. The expressions $a \wedge b := \min(a,b)$ and $a \vee b := \max(a,b)$ stand for the minimum and the maximum  functions respectively. Letters $p,q$ always denote primes; $\P$ denotes the set of all primes. We write $[n,m]=\text{lcm}(n,m)$ and $(n,m) = \text{gcd}(n,m)$ for the least common multiple and greatest common divisor of $m$ and $n$ respectively.

For $b \in \N$ function $\tau_b(n) = \sum_{d_1d_2\cdots d_b= n}{ 1}$ is the $b$-fold divisor function; the standard divisor function is denoted by $\tau(n) = \sum_{d \mid n}{ 1}$. Functions $\varphi(n)$ and $\mu(n)$ are the Euler's totient and the M\"{o}bius function respectively. We use the following functions counting prime factors of $n$:
\begin{equation}
    \Omega(n) := \sum_{\substack{p^\ell \| n}}{\ell},\quad \Omega^*(n) := \sum_{\substack{p^\ell \| n\\ p > 2}}{\ell},\quad \Omega(n,t) := \sum_{\substack{p^\ell \| n \\ p \leq t}}{\ell}, \quad \Omega^*(n,t) := \sum_{\substack{p^\ell \| n \\ 2 <p \leq t}}{\ell}\,.
\end{equation}

Let $b \in \N$ and $\eta >0$. Denote $\mathcal{F}_{b,\eta}$ the class of sub-multiplicative functions satisfying
\begin{equation}\label{def:F b eta}
    \mathcal{F}_{b,\eta} = \{f:\N \to [0,+\infty): 0 \leq f(n) \leq \tau_b(n) \eta^{\Omega(n)},\, f(nm) \leq f(n)f(m)\ \forall n,m \in \N \}\,.
\end{equation}

For $\rho > 0$ we define $ Q(\rho) := \rho \log (\rho/e) + 1$
such that $\gamma_0:= Q(2) = \log 4 - 1 \approx 0.386$ and
\begin{equation}
    \delta := Q(1/\log 2) = 1 - \frac{1+\log\log 2}{\log 2}\approx 0.086\,.
\end{equation}

Let $3 \le y < z \le x$, and let $z^\prime = \min\{z,x/y\}$. Throughout the paper for an interval $(y,z]$ parameters $u,\gamma$ are always defined by
\begin{equation}
    u := u(y,z) = \frac{\log (z/y)}{\log z},\quad \gamma := \gamma(x,y,z) =\frac{\log (1/u)}{\log\log z^\prime}\,.
\end{equation}

Finally, we denote the number of ways to write $n$ as a product of an integer and a shifted prime from $(y,z]$ by
\begin{equation}\label{defn:rep_fn}
    r(n) := r(n,y,z) = \#\{(m,p): n = m(p-1),\ p \in (y,z] \} = \sum_{p\in (y,z]} {\bm{1}_{(p-1)\mid n}}\,,
\end{equation}
and the number of integers $n \le x$ with at least one such representation by
\begin{equation}
    H^*(x,y,z) := \#\{n \le x: \exists p \in (y,z],\ p-1 \mid n\} = \sum_{n \leq x}{\bm{1}_{r(n) \ge 1}}\,.
\end{equation}

%% file: 02_technical_details.tex
For our proof we largely adopt the approach of \cite{tenenbaum1984probabilite} and 
\cite{ford2017integers}.

One can bound $H^*(x,y,z)$ in terms of its first and second moments $M_1$ and $M_2$ defined by 
\begin{equation}\label{defn:rep_fn_moments}
        M_1 := \sum_{n \leq x} r(n), \quad M_2 := \sum_{n \leq x} r(n)^2\,.
\end{equation}

By Markoff's inequality (``first moment method'') we have the upper bound $H^*(x,y,z) \le M_1$. The lower bound $H^*(x,y,z) \ge M_1^2/M_2$ follows from the Paley--Zygmund inequality (``second moment method''). Thus, combining the two we have
\begin{equation}\label{eq:H_M1M2_bound}
    \frac{(M_1)^2}{M_2} \le H^*(x,y,z) \le M_1\,.
\end{equation}
Our proof boils down to utilising these bounds. However, depending on the regime, the set of integers $n$ we consider has to be restricted in order to avoid sparse subsets of the integers with a disproportionately large contribution to $M_1$ and $M_2$.

\subsection*{Part(i): short intervals.} In \S\ref{sec:short}, we treat intervals $(y,z]$ that are relatively short such that
\begin{equation}
    \gamma =  \frac{\log (1/u)}{\log\log z^\prime} > \gamma_0 = \log 4 - 1
\end{equation}
is a fixed constant. In this regime $\bm{1}_{r(n) \ge 1}$ and $r(n)$ almost always behave in the same way so that $H^*(x,y,z) = M_1 (1+ o(1))$. In other words, intervals $(y,z]$ typically contain a single shifted prime divisor if they contain one at all. This behaviour persists even in the vicinity of $\gamma_0$ as long as $\xi = (\gamma-\gamma_0) \sqrt{\log \log z^\prime} \to +\infty$.

We prove the following more precise version of \Cref{mainthm}(i).
\begin{theorem} \label{thm:short intervals}
Let $3 \le y < z \le x$, and let $z^\prime = \min\{z,x/y\}\to \infty$ as $x \to \infty$. Assume that $\xi = (\gamma-\gamma_0)\sqrt{\log\log z^\prime} \to +\infty$, then 
\begin{equation}
     x^{-1}H^*(x,y,z) = \left( \sum_{y < p \le z} \frac{1}{p-1}\right) (1+ o(1))\,.
\end{equation}
Further, for $z$ satisfying $y + y^{17/30+\varepsilon} \le z \le y^{1+\varepsilon}$ we have $x^{-1}H^*(x,y,z) \sim u$.
\end{theorem}

For genuinely short intervals with $z \le y + y^{1-\varepsilon}$ inequality \eqref{eq:H_M1M2_bound} is asymptotically sharp. For such short intervals, one can bound the number of pairs $p-1,q-1 \mid n$ by pairs $p-1,d\mid n$, where only one of the divisors is required to be a shifted prime. Using this crude estimate we show that $M_2 \le M_1(1+o(1))$ and thus $H^*(x,y,z) \ge M_1 (1+ o(1))$ regardless of whether we can compute $M_1$ exactly and whether it is even non-zero. 

When $\log(z/y) \gg 1$ we cannot rely on \eqref{eq:H_M1M2_bound} anymore. One does not need to control the anatomy of $n$ to compute $M_1$ (see Lemma \ref{lem:M1}), but the main contribution to $M_1$ indeed comes from the integers satisfying \eqref{eq:Omeganz_i 0} as shown in \Cref{prop:M1mild}. However, $M_2$ becomes dominated by even more unusual integers with $\Omega(n,z^\prime) \ge b \log\log z^\prime$ for some $b > 2$. 

Thus, for the purposes of obtaining the lower bound we restrict our set of integers accordingly to those satisfying a mild cut-off condition
\begin{equation}\label{def:mild cutoff w}
    \Omega(n,z^\prime) \le w := 2 \log \log z^\prime + \psi(z^\prime) \sqrt{\log\log z^\prime}
\end{equation}
for some slowly growing function $\psi(t)$ such that $\psi(t) \le \min\{(\log\log t)^{1/2-\varepsilon},\xi\}$. Other than that, the choice of $\psi(t)$ does not matter: we only need to make sure that the underlying Poisson random variable with parameter $\lambda^\prime = 2\log\log z^\prime$ coming from the first moment described in \eqref{eq:shifted_prime_heuristic0} stays within the Gaussian regime around its expectation.

In \Cref{prop:M1mild} we confirm that condition \eqref{def:mild cutoff w} does not affect the first moment, and in \Cref{prop:M2 short intervals} we show that this condition prohibits pairs of shifted prime divisors so that
\begin{equation}\label{defn:rep_fn_moments_w}
        \tilde{M}_1 := \sum_{n \leq x} r(n) \bm{1}_{\Omega(n,t)\le w} \sim M_1, \quad \tilde{M}_2 := \sum_{n \leq x} r(n)^2 \bm{1}_{\Omega(n,t)\le w} \sim M_1\,.
\end{equation}
Both of the these facts are proved via Chernoff inequalities. Therefore, we establish that
\begin{equation}\label{eq:H_lbound_w}
 H^*(x,y,z) \ge \frac{\Tilde{M_1}^2}{\Tilde{M_2}} \sim M_1\,.
\end{equation}

 As $\gamma$ approaches $\gamma_0$ pairs of shifted prime divisors become more common. We observe that there is a continuous transition at $\gamma=\gamma_0$ marking a qualitative change in the anatomy of integers producing shifted prime divisors. 

\subsection*{Parts (ii) and (iii): longer intervals}

Next, in \S\ref{sec:long} we consider intervals $(y,z]$ such that 
\begin{equation}
    \varepsilon_1 < \gamma  \le \gamma_0 + \varepsilon_2
\end{equation}
for some small $\varepsilon_1,\varepsilon_2 > 0$. That is, we treat $0<\gamma \le \gamma_0$ together with the transition around $\gamma = \gamma_0$. Conceptually, this is the most challenging part of the paper. 

We prove the following version of \Cref{mainthm}(ii) and \Cref{mainthm}(iii).

\begin{theorem} \label{thm:longer intervals}
Let $3 \leq y < z \leq x$ such that $z^\prime = \min\{z, x/y\} \to \infty$ as $x \to \infty$. Let $\varepsilon > 0$ be small and fixed. Let $k_0 = \frac{1+\gamma}{\log 2} \log \log z^\prime$. Then for $z \ge 2y$ such that $\gamma \ge \varepsilon$, we have
\begin{equation}
        x^{-1}H^*(x,y,z) \asymp \frac{1}{\left(\log{z^\prime}\right)^{ \gamma}}\sum_{k\le k_0}{\frac{\left(2\log\log z^\prime \right)^k}{k!(\log z^\prime)^2}}\,,
\end{equation}
where only the implied constant in the lower bound depends on $\varepsilon$. Further, we have the following.
    \begin{enumerate}[(i)]
    \item If $\xi = (\gamma-\gamma_0)\sqrt{\log\log z^\prime} \ge -C$ for some large positive constant $C$, then 
\begin{equation}
   x^{-1}H^*(x,y,z) \asymp u\,.
   \end{equation}
\item If $\xi = (\gamma-\gamma_0)\sqrt{\log\log z^\prime}\to - \infty$ as $x \to \infty$, then
        \begin{equation}
    \begin{aligned}
        x^{-1}H^*(x,y,z) \asymp
                (\log z^\prime)^{ -Q(\frac{1+\gamma}{\log 2})} ((\gamma_0-\gamma)\sqrt{\log\log z^\prime} + 1)^{-1}\,,
    \end{aligned}
\end{equation}
where only the implied constant in the lower bound depends on $\varepsilon$. 
\end{enumerate}
\end{theorem}

Note that the range in \Cref{thm:longer intervals} overlaps with \Cref{thm:short intervals}, however, due to certain technical assumptions we make in the proof this method fails to obtain an asymptotic and only establishes the correct order of magnitude for $\gamma \ge \gamma_0 + \varepsilon$. When $\xi \le 1$, we do not know if an asymptotic exists at all; likely the answer is the same as one for $H(x,y,z)$.

The proof in this regime is different in two ways. 

First of all, for $\gamma \le \gamma_0$ integers dominating $H^*(x,y,z)$ no longer have $\Omega(n,z^\prime) \approx 2 \log \log z^\prime$. Instead, most of the contribution comes from integers with
\[
\Omega(n,z^\prime) \lesssim \frac{1+ \gamma}{\log 2} \log\log z^\prime  \le 2\log \log z^\prime\,.
\]
This means that we will need to fix $\Omega(n,z^\prime) = k$, and apply moment methods separately to each such set to evaluate
\begin{equation}
    H^*_k (x,y,z) := \#\{n \le x: \Omega(n,z^\prime) = k,\ \exists p \in (y,z],\ p-1 \mid n\}\,.
\end{equation}
 For each $k$ by Markoff's inequality we have 
\begin{equation}
    H^*_k(x,y,z) \le M_{1,k} := \sum_{\substack{n\leq x}}{r(n)\bm{1}_{\Omega(n,z^\prime) = k}}
\end{equation}
as well as the trivial bound $H^*_k(x,y,z) \le \pi_k = \#\{n\le x: \Omega(n,z^\prime) = k\}$. We compute $M_{1,k}$ and thus justify \eqref{eq:shifted_prime_heuristic0} in \Cref{lem: M1 fixed k}. This gives us a new compound first moment upper bound 
\begin{equation}
     H^*(x,y,z) = \sum_{k\ge 0} H^*_k(x,y,z) \le \sum_{k\ge 0} \min \{M_{1,k}, \pi_k\}\,,
\end{equation}
which we will show to have the correct order of magnitude for fixed $\gamma > 0$.

If we were to apply the Paley--Zygmund inequality directly for each $\log\log z^\prime \le k \le 2\log\log z^\prime$, this approach would fail us as 
\begin{equation}
    M_{2,k}:= \sum_{\substack{n\leq x}}r^2(n)\bm{1}_{\Omega(n,z^\prime) = k}
\end{equation}
is dominated by integers with $\Omega(n,t^\prime) \ge 2 \log\log t^\prime$ for some $t^\prime \le z^\prime$. This was acceptable for $\gamma \ge \gamma_0$ when $k \approx 2\log\log z^\prime$, but for $\gamma < \gamma_0$ such integers are rare.

Suppose $\Omega(n,t) = k = \rho \log\log t$ is fixed, then for any $t^\prime \le t$ on average $\Omega(n,t^\prime) \approx \rho \log\log t^\prime$. Further, for most integers with $\Omega(n,t) = k$ 
\begin{equation}\label{eq:barrier_observation}
    \Omega(n,t^\prime) - \rho \log\log t^\prime  \ll_\varepsilon  \min\{\log\log t^\prime, \log\log t - \log\log t^\prime\}^{1/2+\varepsilon}
\end{equation}
for any $t^\prime \le t$ and $\varepsilon > 0$. It is then natural to restrict to integers satisfying \eqref{eq:barrier_observation}. 

Let $(z^\prime)^c \le t \le z^\prime$. We define the following barrier-type condition for $\{\Omega(n,t^\prime)\}_{t^\prime \le t}$. Let $t_j = e^{e^j}$ for $0 \le j \le \log\log t$. Let $g(s) = s^a$ for some $1/2+\varepsilon \le a \le 1-\varepsilon$. Define the barrier
\begin{equation} \label{def:barrier outline 0}
    B(\ell) = \rho (\ell + g(\min\{\ell, \log\log t - \ell\}))\,
\end{equation}
and restrict our set of integers to $n$ satisfying the barrier condition
\begin{equation}\label{def:barrier outline}
    \mB = \mB(k) = \{n \in \N: \Omega(n,t)= k,\ \Omega(n, t_j) \le B(j) +C\ \forall j \in [0,\log\log t]\}\,.
\end{equation}
While the set of exceptional $n$ failing the barrier condition is not necessarily vanishing (unless we choose $C$ to grow with $t$), by choosing $g(s) \ge s^{1/2+\varepsilon}$ and $C$ sufficiently large, we restrict to a positive proportion of integers within the set where $\Omega(n,t)=k$.

Applying a separate barrier for each $k \le k_0$, we replace $H^*_k(x,y,z)$ with
\begin{equation}
    \tilde{H}^*_k(x,y,z) = \sum_{n \le x} {\bm{1}_{r(n) \ge 1} \bm{1}_{n \in \mB(k)}} \ge \frac{M_{1,k}(\mB)^2}{M_{2,k}(\mB)}\,,
\end{equation}
where $M_{1,k}(\mB)$ and $M_{2,k}(\mB)$ are the moments under the barrier condition.

In \Cref{prop:firstmomentbarrier} we show that \eqref{def:barrier outline} indeed does not affect the order of magnitude of the first moment such that $M_{1,k}(\mB) \asymp M_{1,k}$. This is done via a union bound together with Chernoff's inequality for each $j\le \log\log t$. 

The key proposition for this regime is \Cref{prop:second moment general barrier} where we show that for an arbitrary barrier, not necessarily defined as above, the second moment exhibits a factor proportional to the sum of all partial first moments. This bound is essentially sharp for ``reasonable'' barriers such as \eqref{def:barrier outline 0}.

For our choice \eqref{def:barrier outline 0} \Cref{prop:second moment general barrier} implies that for $k \le \min\{2,\frac{1+\gamma}{\log 2}\}\log\log t$
\begin{equation}
    M_{2,k}(\mB) - M_{1,k}(\mB) \ll M_{1,k} \cdot \frac{\log (z/y)}{\log z} \sum_{\log\log (z/y) \le j\le \log\log t} 2^{B(j)} e^{-j}\,.
\end{equation}
The final term in the sum above is of order $2^k/\log z^\prime$, so this bound is at least of order $M_{1,k}^2/\pi_k(x)$, which is of course consistent with H\"{o}lder's inequality. Thus, if the sum is dominated by the final term, the bound above is in fact sharp. In particular, this occurs for fixed $\gamma > 0$ since
\begin{equation}
    B(j)\log 2 - j = (\rho \log 2 - 1) j + \rho g(j) \le \gamma j + 2 g(\min\{j,\log\log t - j\})
\end{equation}
is monotone increasing. Hence for fixed $\gamma > 0$ the lower bound matches the upper bound for all $k_0(1-\varepsilon) \le k \le k_0$, which then in turn holds for their respective sums.

\subsection*{Part (iv): hidden transition} 
For $0 \le \gamma \le \gamma_0- \varepsilon$ the lower and the upper bounds for $H^*(x,y,z)$ are dominated by $H^*_k(x,y,z)$ with $k = \frac{1+\gamma}{\log 2} \log \log t$. This implies that the lower bound matches the upper bound up to a factor
\begin{equation}
    V = \frac{\log (z/y)}{\log z} \sum_{\log\log (z/y) \le j\le \log\log t} e^{\gamma j + \rho g(j)} = \sum_{j\le \lambda_0} e^{-\gamma j + \rho g(j)}\,.
\end{equation}
If $\gamma \gg 1$, the sum above is convergent, thus $V \ll 1$. If $\log (z/y) \asymp \log (x/y)$, then $\lambda_0 \le \log \log (x/y) - \log \log (z/y) \ll 1$, so this sum only has finitely many terms, all of order $O(1)$. This includes the case $u \gg 1$ and therefore $\gamma \gg 1/\log\log z^\prime$.

In the remaining range where $\gamma =o(1)$ such that $\gamma \gg 1/\log\log z^\prime$ and $\lambda_0 \to \infty$, our lower bound contains a non-trivial additional factor. 

With this in mind, we state the following version of \Cref{mainthm}(iv).

\begin{theorem} \label{thm:gamma to 0}
Let $3 \leq y < z \leq x$ such that $z^\prime = \min\{z, x/y\} \to \infty$ as $x \to \infty$. Let $\lambda_0 = \log\log z^\prime - \log\log (z/y)$. Then for $0\le \gamma \le \gamma_0-0.01$ and any fixed $\varepsilon \in (0,1/2)$ we have
\begin{equation}
        H^*(x,y,z) \gg_\varepsilon \frac{x}{(\log z^\prime)^{Q(\frac{1+\gamma}{\log 2})} \sqrt{\log\log z^\prime}} \cdot V(\varepsilon)^{-1}\,,
\end{equation}
where
\begin{equation}
    V(\varepsilon) := \sum_{j\le \lambda_0} e^{-\gamma j + j^{1/2+\varepsilon}} \ll \exp\left(C(\varepsilon) (\min\{\lambda_0 + 1,1/\gamma^2\})^{1/2+2\varepsilon} \right)\
\end{equation}
for some positive constant $C(\varepsilon)>0$ depending on $\varepsilon$.
\end{theorem}
We prove \Cref{thm:gamma to 0} together with \Cref{thm:longer intervals} in \S\ref{sec:long}.

%% file: 03_preliminaries.tex
In this section we present a few standard estimates for sums of multiplicative functions and Poisson statistics as well as some tools from sieve theory. Most of the results in this section are well established in the literature or happen to be standard exercises.

\subsection{Sums of multiplicative functions and tools from sieve theory}

The following lemma is a standard estimate on sums of multiplicative functions, see for example {\cite[Théorème 3.5]{tenenbaum2022introduction}}.

\begin{lemma}\label{lem:Tenenbaumsum}
    Let $f \in \mathcal{F}_{b,\eta}$ for some $b \in \N$ and $\eta \ge 1$. Then, for all $x\ge 2$, we have
    \[ 
    \sum_{n\le x} f(n)\ll_{b,\eta} \frac{x}{\log{x}}\sum_{n\le x}\frac{f(n)}{n}\,. 
    \]
Further, if $\eta \in [1,2)$, then
\begin{equation}
    \sum_{n \le x} \frac{f(n)}{n} \ll_{b,\eta} \exp \left(\sum_{p \le x}\frac{f(p)}{p}\right)\,. 
\end{equation}
\end{lemma}

The following lemma is a classical sieve estimate; see \cite[Theorem 7.11]{montgomery2007multiplicative} or \cite[Exercise 14.5]{koukoulopoulos2020distribution} for an alternative proof method.

\begin{lemma}
\label{lem:rough}
    Uniformly for $3/2 \le t \le x/2$, we have $\#\{ n \leq x : P^-(n)>t \} \asymp \frac{x}{\log{t}}$. The upper bound holds for $3/2 \le t \le x$.
\end{lemma}

The following lemma is a variant of Lemma \ref{lem:Tenenbaumsum} for short intervals adapted to shifted primes. It follows directly from {\cite[Lemma 2.5]{koukoulopoulos2010divisors}}.

\begin{lemma}
\label{lem:Koukoulopoulossum}
    Let $f \in \mFb$ for some $b \in \N$ and $\eta \ge 1$. Assume additionally that there exists a collection $\mathbf{B} = \{B_{f,\varepsilon}\}_{\varepsilon>0}$ of positive real numbers such that for all $m,n \geq 1$
    \begin{equation}
        f(mn) \leq \left(\min_{\varepsilon > 0} B_{f,\varepsilon}m^{\varepsilon}\right) f(n)\,.
    \end{equation}
    Then, for any integer $d\le x^{1-\varepsilon}$ and $h$ such that $d x ^\varepsilon  \le h \le x$
    \[
        \sum_{\substack{x-h < p \leq x \\ p \equiv 1 (d)}} f\left( \frac{p-1}{d} \right)
        \ll_{\eta,\mathbf{B}} \frac{h}{\varphi(d)(\log{x})^2}  \sum_{n \leq x} \frac{f(n)}{n}\,.
    \]
\end{lemma}

The following lemma is a generalization of \cite[Lemma 4.4]{ford2017integers} (where $f \equiv 1$). It is proved via a standard application of Brun's sieve as outlined in {\cite[Lemma 4.4]{ford2017integers}}. Alternatively, one can derive the statement from \cite[Lemma 4.4]{ford2017integers} by taking out $P^+(n)$, the largest prime factor of $n$. 

\begin{lemma}
\label{lem:Fordsieve} Let $f \in \mathcal{F}_{b,\eta}$ for some $b \in \N$ and $\eta \in [1,2)$. Let $A<B$ be even and positive integers. Then uniformly for all $3/2 \le t \le x/2$ we have
% \begin{equation}
%     \sum_{\substack{n \le x\\ P^-(n)>t\\An+1 \in \mathbb{P}}} f(n) \ll_{b,\eta} \frac{A}{\varphi(A)} \cdot \frac{x}{(\log x)^2} \exp\left(\sum_{t < p\le x} \frac{f(p)}{p}\right)
% \end{equation}
\begin{equation}
    \sum_{\substack{n \le x\\An+1 \in \mathbb{P}}} f(n) \ll_{b,\eta} \frac{A}{\varphi(A)} \cdot \frac{x}{(\log x)^2} \exp\left(\sum_{p\le x} \frac{f(p)}{p}\right)
\end{equation}
    and 
% \begin{equation}
%     \sum_{\substack{n \le x\\ P^-(n)>t\\An+1,Bn+1 \in \mathbb{P}}} f(n) \ll_{b,\eta} \Phi(A,B) \cdot \frac{x}{(\log x)^3} \exp\left(\sum_{t < p\le x} \frac{f(p)}{p}\right)
% \end{equation}
\begin{equation}
    \sum_{\substack{n \le x\\An+1,Bn+1 \in \mathbb{P}}} f(n) \ll_{b,\eta} \Phi(A,B) \cdot \frac{x}{(\log x)^3} \exp\left(\sum_{p\le x} \frac{f(p)}{p}\right)
\end{equation}
    where
    \begin{equation}\label{lem:Fordsieve-0}
        \Phi(A,B):=\prod_{p|AB(B-A)} \frac{p}{p-1} \prod_{p|(A,B)} \frac{p}{p-1} = \frac{AB(B-A)}{\varphi(AB(B-A))}\cdot\frac{(A,B)}{\varphi(A,B)}.
    \end{equation}
\end{lemma}

The first part of \Cref{lem:Fordsieve} partly overlaps with \Cref{lem:Koukoulopoulossum}. \Cref{lem:Koukoulopoulossum} relies on existing estimates for primes in arithmetic progressions on short intervals, so the class of admissible functions $f$ is more restricted, in particular, it does not admit $f(n) = \eta^{\Omega(n,t)}$ for $\eta > 1$. This restriction is only relevant on short intervals, so it does not come up in \Cref{lem:Fordsieve}. We will be using \Cref{lem:Koukoulopoulossum} only in \S\ref{sec:short} when dealing with $z \le 2y$. In all other cases we will be using \Cref{lem:Fordsieve} instead due to it being more flexible.

\subsection{Poisson statistics}
We will require estimates on the number of integers with fixed number of prime factors. The distribution of $\Omega(n)$ computed by Sathe \cite{sathe1953problem} and Selberg \cite{selberg1954note} is known to roughly follow the Poisson distribution with parameter $\lambda = \log\log n$. Further, Hal\'asz \cite{halasz1972remarks} showed that a similar statement holds for $\Omega(n,t)$ and $\Omega^*(n,t)$.

\begin{lemma}[{\cite{halasz1972remarks}}]\label{lem:Big-Omega-bounds}
Let $\varepsilon>0$ and $3/2 < t < x$. Then for all integers $0 \le k\leq (2-\varepsilon)\log\log t$ we have
\[
   x^{-1} \#\{n \leq x :\Omega(n,t)=k\} \ll_{\varepsilon}\frac{(\log\log t)^k}{k!\log t}\,.
\]
If $\varepsilon\log\log t<k\leq (2-\varepsilon)\log\log t$, we also have
\[
x^{-1}\#\{n \leq x :\Omega(n,t)\in \{k,k+1\}\} \gg_{\varepsilon}\frac{(\log\log t)^k}{k!\log t}\,.
\]
\end{lemma}
Note that corresponding statements also hold in logarithmic density as follows from their respective proofs.

It is a standard exercise to show a similar Poisson-type upper bound for sums of non-negative multiplicative functions. We provide a sketch of the proof below based on the approach of Hall and Tenenbaum \cite[Theorem 08]{Hall_Tenenbaum_1988}.

\begin{lemma} \label{lem: hardy-ramanujan multiplicative} 
Let $f \in \mathcal{F}_{b,\eta}$ for some $b \in \N$ and $\eta \in [1,2)$, and let $\varepsilon>0$ be small. Let $\lambda(f) = \sum_{p \leq t} f(p)/p$. Then uniformly over all $t \in (2,x]$ and integers $1 \le k \le \frac{2- \varepsilon}{\eta} \lambda(f)$
one has  
      \[
            S_f(x,t;k)= \frac{1}{x} \sum_{\substack{n\le x \\ \Omega(n,t) = k}} f(n) \ll_{b,\eta,\varepsilon} \frac{1}{k!\log x} \brackets{\sum_{p \leq t}\frac{f(p)}{p}}^k \exp \bigg(\sum_{t <p \leq x}\frac{f(p)}{p}\bigg).
   \] 
\end{lemma}
\begin{proof}
We start by moving to logarithmic density. Write $\log x = \sum_{p\mid n}\log p + g(n)$, where $g(n) := \log x - \sum_{p\mid n}\log p$, then
\begin{equation}
    \log x \sum_{\substack{n\le x \\ \Omega(n,t) = k}} f(n) = \sum_{\substack{n\le x \\ \Omega(n,t) = k}} f(n) \sum_{p | n}\log p + \sum_{\substack{n\le x \\ \Omega(n,t) = k}} f(n) g(n) =: x(S_1+S_2)\,.
\end{equation}
Let us bound $S_1$. Since $f \in \mathcal{F}_{b,\eta}$, we have $f(p) \leq 2b \eta $, so
\begin{equation}\label{eq proof l2.7 s1}
\begin{aligned}
        S_1 & \le \frac{1}{x}\sum_{p \leq x}f(p) \log p \sum_{\substack{m \leq x/p \\ \Omega(mp,t) = k}}f(m) \leq \frac{1}{x} \sum_{\substack{m \leq x \\ \Omega(m,t) \in \{k-1,k\}}}f(m) \sum_{p \leq x/m}f(p) \log p\\
        & \ll_{b,\eta} \sum_{\substack{m \leq x \\ \Omega(m,t) = k-1}}\frac{f(m)}{m}+\sum_{\substack{m \leq x \\ \Omega(m,t) = k}}\frac{f(m)}{m} \,.
\end{aligned}
\end{equation}
Let $\lambda(f) = \sum_{p \leq t} f(p)/p$. Given a complex number $z$ satisfying $|z| \leq (2-\epsilon)/\eta$, we have
\begin{equation} \label{eq proof l2.7 3}
    \begin{aligned}
        \sum_{P^+(m) \leq x}\frac{f(m)z^{\Omega(m,t)}}{m} &= \prod_{p \leq t}\bigg(1+\frac{z f(p)}{p}+\frac{z^2 f(p^2)}{p^2}+\cdots\bigg)\prod_{t < p \leq x}\bigg(1+\frac{f(p)}{p}+\cdots\bigg) \\ &\ll_{b,\eta} \exp{\bigg(\Re(z) \lambda(f)\bigg)}\exp{ \bigg(\sum_{t< p \leq x}\frac{f(p)}{p}\bigg)}
    \end{aligned}
\end{equation}
We write $z = r e^{i \theta}$ where $r = k/\lambda(f)$. By Cauchy's theorem and \eqref{eq proof l2.7 3}, one has
\begin{equation}
    \begin{aligned}
           \sum_{\substack{m \leq x \\ \Omega(m,t) = k}}\frac{f(m)}{m} &\le \frac{1}{2 \pi r^k} \int_{0}^{2\pi } e^{-i k \theta}\sum_{P^+(m) \leq x}\frac{f(m)z^{\Omega(m,t)}}{m} \;d\theta\\
   &\ll_{b,\eta} \frac{1}{2 \pi r^k}  \exp{ \bigg(\sum_{t< p \leq x}\frac{f(p)}{p}\bigg)}\int_{0}^{2\pi} e^{ k \cos \theta}\;d\theta\\
   &\ll_{b,\eta} \frac{\lambda^k}{k^k} \exp{ \bigg(\sum_{t< p \leq x}\frac{f(p)}{p}\bigg)}\int_{0}^{2\pi} e^{k \cos \theta}\;d\theta.
    \end{aligned}
\end{equation}
The integral in the last line of the above equation is $\ll k^{-1/2}e^k,$ and so by Stirling's formula 
\begin{equation} 
    \sum_{\substack{m \leq x \\ \Omega(m,t) = k}}\frac{f(m)}{m} \ll_{b,\eta} \frac{\lambda(f)^k}{k!}\exp{ \bigg(\sum_{t< p \leq x}\frac{f(p)}{p}\bigg)} \label{eq proof 12.7 log average upper bound}\,.
\end{equation}
Returning to \eqref{eq proof l2.7 s1} and noting that a similar bound applies to the sum with $\Omega(n,t) = k-1$, we conclude that $S_1$ satisfies the required bound
\begin{equation} 
    S_1  \ll_{b,\eta}  \frac{\lambda(f)^k}{k!\log x} \bigg(1 + \frac{k}{\lambda(f)}\bigg) \exp{ \bigg(\sum_{t< p \leq x}\frac{f(p)}{p}\bigg)}  \ll_{b,\eta} \frac{\lambda(f)^k}{k!\log x} \exp{ \bigg(\sum_{t< p \leq x}\frac{f(p)}{p}\bigg)}\,.
\end{equation}

The last step in the proof is to show that $S_2$ is sufficiently small compared to $S_1$. The argument of Hall and Tenenbaum \cite[p.6]{Hall_Tenenbaum_1988} shows that when $f \equiv 1$ one has
\begin{equation}
   \frac{S_2}{\log x} \ll \frac{S_f}{\sqrt{\log x}}
\end{equation}
using the fact that $\sum_{n \le x} g(n)^m \ll_m x$  for any $m \in \N$. In our case, the upper bound for $S_f(k)$ in the statement of the Lemma is at least of order $1/\log x$, so we can assume that $S_f(k) \gg 1/\log x$, otherwise the required bound is trivial.

Recall that $f(p) \le 2b \eta$.  Let $a \in (\frac{\log \eta}{\log 2},1)$, then by \Cref{lem:Tenenbaumsum} we have
\begin{equation}
    \begin{aligned}
        \frac{1}{x} \sum_{n \le x} f(n)^{1/a} \ll_{a,b,\eta} \exp\brackets{\sum_{p \le x} \frac{f(p)^{1/a}-1}{p}} \ll (\log x)^{(2b \eta)^{1/a}-1}\,.
    \end{aligned}
\end{equation}
Let $\mu \in (0,1)$, then by Hölder's inequality for exponents $a_1 = 1/(1-\mu)$, $a_2 = 1/(\mu a)$ and $a_3 = 1/(\mu (1-a))$,
\begin{equation}
    \begin{aligned}
   \frac{S_2}{\log x}  & = \frac{1}{x\log x} \sum_{n \le x}f(n)g(n) \bm{1}_{\Omega(n,t) = k} \\
        & \le \frac{1}{x\log x}  \brackets{\sum_{n \le x} f(n)\bm{1}_{\Omega(n,t) = k}}^{1-\mu} \brackets{\sum_{n \le x} f(n)^{1/a}}^{\mu a} \brackets{\sum_{n \le x} g(n)^{1/(\mu(1-a))}}^{\mu(1-a)} \\
        & \ll_{\mu,a} S_f(k) \brackets{\log x}^{\mu + \mu a ((2b \eta)^{1/a}-1)-1} \,.
    \end{aligned}
\end{equation}
Take $\mu$ such that $\mu (1 + a((2b \eta)^{1/a}-1)) < 1$, and the statement of the lemma follows.
\end{proof}

Shifted primes---more precisely, the numbers $(p-1)/2$ where $p$ is prime, are expected to have similar arithmetic properties to generic integers. In particular, their respective prime factors follow roughly the same distribution up to a factor corresponding to the density of primes as shown by Timofeev \cite{timofeev1995hardy}. The upper bound is proved in \cite[Theorem 2]{timofeev1995hardy}. The matching lower bound is a special case of \cite[Theorem 3]{timofeev1995hardy}, however as pointed out in \cite{ford2017integers}, there is a superfluous assumption in \cite{timofeev1995hardy} which must be removed. See the proof of \cite[Lemma 3.5]{ford2017integers} for more details.

\begin{lemma}[{\cite{timofeev1995hardy}}]
\label{lem:Big-Omega-bounds-shifted}
Let $\varepsilon>0$ and $3/2 < t < x$. Then, for all integers $0\le k\leq (2-\varepsilon)\log\log t$, we have
\[
    \left(\frac{x}{\log x}\right)^{-1}\#\{p \leq x: \Omega(p-1,t)=k\} \ll_{\varepsilon}\frac{(\log\log t)^k}{k!\log t}\,.
\]
Further, there exists a constant $c = c(\varepsilon) >1$ such that for all real numbers $c(\varepsilon) < t < x$ and integers $\varepsilon \log\log t \le k\leq (2-\varepsilon)\log\log t$, we have
\begin{equation}
        \left(\frac{x}{\log x}\right)^{-1} \#\{p \leq x: \Omega(p-1,t)\in \{k,k+1,k+2\}\} \gg_{\varepsilon} \frac{(\log\log t)^k}{k!\log t}\,.
\end{equation}
\end{lemma}

A similar upper bound holds for sums over shifted primes weighted by a non-negative multiplicative function.

\begin{lemma} \label{lem: hardy-ramanujan multiplicative shifted prime} Let $f \in \mathcal{F}_{b,\eta}$ for some $b \in \N$ and $\eta \in [1,2)$, and let $\varepsilon>0$ be small. Let $\lambda(f) = \sum_{p \leq t} f(p)/p$. Then uniformly over all $t \in (2,x]$ and $0 \le k \le  \frac{2- \varepsilon}{\eta} \lambda(f)$ one has 
      \[
            S_f^*(x,t;k) =  
            \left(\frac{x}{\log x}\right)^{-1} \sum_{\substack{p\le x \\ \Omega(n,t) = k}} f(p-1) \ll_{b,\eta,\varepsilon} \frac{1}{k!\log x} \brackets{\sum_{p \leq t}\frac{f(p)}{p}}^k  \exp \bigg(\sum_{t <p \leq x}\frac{f(p)}{p}\bigg)\,.
   \] 
\end{lemma}
\begin{proof}
% Let $\rho = k/\lambda(f)\le (2-\varepsilon)/\eta$.

Let $q = P^+(p-1)$ be the largest prime factor of $p-1$, then $p-1 = d q$, where $P^+(d) \le q$. Let $P_0 \le t$ be the largest integer such that
\begin{equation}\label{eq: P0}
    \frac{2-\varepsilon/2}{\eta}\sum_{p \le P_0} \frac{f(p)}{p} \le k \le \frac{2-\varepsilon}{\eta}\lambda(f)\,.
\end{equation}
Write
\begin{equation}\label{eq:HR mult shifted 1}
     \frac{x}{\log x} \cdot S_f^*(x,t;k) \le \sum_{\substack{\sqrt{x} \le p\le x\\ P^+(p-1) > P_0\\ \Omega(d,t) =k}} f(p-1) + \sum_{\substack{\sqrt{x} \le p\le x \\ P^+(p-1)\le P_0}} f(p-1) + \sum_{\substack{p\le \sqrt{x}}} f(p-1) \,.
\end{equation}
By \Cref{lem:Fordsieve} the contribution of $p\le \sqrt{x}$ is $O(x^{1/2+\varepsilon})$, so from now assume that $p \ge \sqrt{x}$. 

Denote the sum over $p$ with $P^+(p-1) > P_0$ in \eqref{eq:HR mult shifted 1} by $S_1$, and over $p$ with $P^+(p-1) \le P_0$ by $S_2$. We claim that the contribution of $p$ with $P^+(p-1)\le P_0$ is negligible. 

Let $\xi_0 = 1/(10j_0)$. By \Cref{lem:Fordsieve}
\begin{equation}
    \begin{aligned}
        S_2 & \le 
        x^{-\xi_0} \sum_{\substack{\sqrt{x} \le p\le x \\ P^+(p-1)\le P_0}} f(p-1) (p-1)^{2\xi_0}  \ll_{b,\eta} \frac{x^{1-\xi_0} }{\log^2 x} \exp \brackets{\sum_{p\le P_0}\frac{f(p)}{p^{1-2\xi_0}}} \ll \frac{x^{1-\xi_0}e^{\lambda(f)}}{\log^2 x}\,.
    \end{aligned}
\end{equation}
The upper bound for $S_f^*$ in the statement of the Lemma is at least of order $1/\log x $, so it suffices to check that $x^{-\xi_0}e^{\lambda(f)} = o(1)$ or $\log P_0 = o(\log x/\lambda(f))$. 

Since $f \in \mathcal{F}_{b,\eta}$ we have $f(p) \le 2b \eta$ for any $p$. Together with \eqref{eq: P0} this gives us
\begin{equation}  
\begin{aligned}
        \frac{\varepsilon}{4} \cdot \lambda(f)\le \sum_{P_0 < p \le t} \frac{f(p)}{p} \le 2b \eta (\log\log t - \log\log P_0) + O(b)\,.
    \end{aligned}
\end{equation}
It follows that $\log\log P_0\le \log\log t - \frac{\varepsilon}{8b\eta} \lambda(f) + O(1)$ and thus
\begin{equation}\label{eq: HR shifted 1}
        \log P_0 \ll \log t \cdot \exp \left(-\frac{\varepsilon}{16b}\lambda(f)\right) \ll \log x \cdot \exp \left(-\frac{\varepsilon}{16b}\lambda(f)\right)\,,
\end{equation}
which is certainly small enough. 

It remains to bound $S_1$. 

Fix the size of $q =P^+(p-1)$ such that $q\in I_j = (e^{j-1},e^{j}]$ and $P_0 \in I_{j_0}$. Then by \Cref{lem:Fordsieve}
\begin{equation}
\begin{aligned}
         S_1 \le    \sum_{j_0 \le j \le \log x} \sum_{\substack{\sqrt{x} \le p\le x\\ P^+(p-1) \in I_j\\ \Omega(d,t) =k}} f(p-1)  & \le  \sum_{j_0 \le j \le \log x}\sum_{\substack{\sqrt{x} e^{-j} \le d\le 4xe^{-j}\\ P^+(d) \le e^{j}\\ \Omega(d,t) \in \{k,k-1\}}} f(d) \sum_{\substack{q \in I_j\\dq+1 \in \P}} f(q) \\
     & \ll  \sum_{j_0 \le j \le \log x}\frac{e^j}{j^2} \sum_{\substack{\sqrt{x} e^{-j} \le d\le 4xe^{-j}\\ P^+(d) \le e^{j}\\ \Omega(d,t) \in \{k,k-1\}}} f(d) \cdot \frac{d}{\varphi(d)}\,.
\end{aligned}
\end{equation}

Let $g(d) = f(d) d/\varphi(d)$. By definition $g(p^\ell) = f(p^\ell) \cdot \frac{p}{p-1}$, so $g \in \mathcal{F}_{b,\eta}$. Let $t_j = \min\{t,e^j\}$. For $j \ge j_0$ we always have $\sum_{p\le t_j} g(p)/p = \sum_{p\le t_j} f(p)/p + O(b)\le (2-\varepsilon)k$. Thus, we may apply \Cref{lem: hardy-ramanujan multiplicative} to bound the sum above.

Let $\xi = 1/(10j)$ such that $1 \le (de^j/\sqrt{x})^{2\xi} \ll (d/\sqrt{x})^{2\xi}$. By \Cref{lem: hardy-ramanujan multiplicative} for $j_0 \le j \le \log x$
\begin{equation}\label{eq:hardy-ramanujan multiplicative shifted prime 1}
    \begin{aligned}
        x^{-\xi} \sum_{\substack{d\le 4xe^{-j}\\ P^+(d) \le e^{j}\\ \Omega(d,t_j)=k}} g(d) d^{2\xi} 
        & \ll \frac{x^{1-\xi}e^{-j}}{k!\log x} \brackets{\sum_{p\le t_j}\frac{g(p)}{p^{1-2\xi}}}^k \exp\brackets{\sum_{t_j < p\le x}\frac{g(p)}{p^{1-2\xi}}}\,.
    \end{aligned}
\end{equation}
We may replace $g(p)p^{2\xi}$ with $f(p)$ at the total cost of $O(b)$, which is absorbed into the implied constant. Further, observe that if $t_j < t$,
\begin{equation}
    \exp\left(\sum_{t_j < p\le t}\frac{f(p)}{p}\right) \le \left(\frac{\log t}{j}\right)^{2b\eta}\,.
\end{equation}
Summing \eqref{eq:hardy-ramanujan multiplicative shifted prime 1} over $j$ and using
\begin{equation}
    \sum_{\log t\le j \le \log x} \frac{x^{-1/(10j)}}{j^2} + \sum_{j \le \log t} \frac{x^{-1/(10j)}}{j^2} \left(\frac{\log t}{j}\right)^{2b\eta}\ll \frac{1}{\log x}
\end{equation}
we obtain the desired upper bound for $S_1$ and therefore $S_f^*$.
\end{proof}

In order to access precise estimates involving Poisson distribution, it will be useful to control the left and right tails of that distribution. This is known as Bennett's inequality in some settings, which leads to the following improvement in the case of the Poisson distribution.
% See Tao’s blog “An improvement to Bennett’s inequality for the Poisson distribution”, Proposition 3

\begin{lemma}
\label{lem:improvedBenett} Let $Q(\rho) = \rho \log \rho/e + 1$. Let $\lambda>0$. If $-\frac{C}{\sqrt\lambda}\le u<1$ for some constant $C>0$, then
    $$ \sum_{k\le\lambda\left(1-u\right)}\frac{\lambda^k}{k!e^{\lambda}}\ll_C \frac{\exp{\left(-\lambda Q(1-u)\right)}}{\sqrt{\lambda u^2(1-u)+1}}. $$
    If $u \ge 0$, then
    $$ \sum_{k>\lambda(1+u)}\frac{\lambda^k}{k!e^{\lambda}}\ll \frac{\exp{\left(-\lambda Q(1-u)\right)}}{\sqrt{\lambda\min{\{u,u^2\}}+1}}. $$
\end{lemma}
\begin{proof}
Let $f(k) = \lambda Q(k/\lambda)$. By Stirling's formula it suffices to bound corresponding sums of $e^{-\lambda Q(k/\lambda)}/\sqrt{k} = e^{-f(k)}/\sqrt{k}$. The tail estimates follow by taking the minimum between the sum of the geometric series and the estimate from the stationary point method.

 Observe that $f^\prime(\lambda(1\pm u)) = Q^\prime(1\pm u) = \log (1 \pm u)$ and $f^{\prime\prime} (\lambda(1\pm u)) = (\lambda(1\pm u))^{-1}$.  For $u < 1$ we have $|\log (1 \pm u)|  \asymp |u|$, so
    \begin{equation}
    \begin{aligned}
           e^{\lambda Q(1-u)} \sum_{k\le\lambda\left(1-u\right)} \frac{e^{-f(k)}}{\sqrt{k}} & \ll  \frac{1}{\sqrt{(1-u)\lambda}}\min\left\{\frac{1}{|f^\prime(\lambda(1- u))|} + 1,\frac{1}{|f^{\prime\prime} (\lambda(1- u))|^{1/2}}\right\} \\
          & \ll \min\left\{\frac{1}{|u|\sqrt{(1-u)\lambda}},1\right\} \ll \frac{1}{\sqrt{\lambda u^2(1-u)+1}}\,,
                   \end{aligned}
    \end{equation}
and for $u \ge 0$
    \begin{equation}
    \begin{aligned}
           e^{\lambda Q(1+u)} \sum_{k > \lambda\left(1+u\right)} \frac{e^{-f(k)}}{\sqrt{k}} & \ll  \frac{1}{\sqrt{(1+u)\lambda}}\min\left\{\frac{1}{|f^\prime(\lambda(1+ u))|} + 1,\frac{1}{|f^{\prime\prime} (\lambda(1+ u))|^{1/2}}\right\} \\
          & \ll \min\left\{\frac{1}{\min\{u,1\}\sqrt{(1+u)\lambda}},1\right\} \ll \frac{1}{\sqrt{\lambda \min\{u,u^2\}+1}}\,
                   \end{aligned}
    \end{equation}
    as required.
\end{proof}

%% file: 04_short_intervals.tex
In this section we establish an asymptotic for $H^*(x,y,z)$ when $\gamma = \frac{\log (1/u)}{\log\log z^\prime} > \gamma_0$ and $\xi = (\gamma-\gamma_0)\sqrt{\log\log z^\prime} \to +\infty$. This corresponds to \Cref{mainthm}(i), where $\gamma-\gamma_0>0$ is small and fixed, as well as a part of the transition regime, where the asymptotic persists.

\begin{proof}[Proof of \Cref{thm:short intervals}]
We wish to prove that $H^*(x,y,z) \sim M_1$, where by \Cref{lem:M1} 
\begin{equation}
    M_1 := \sum_{n \leq x} r(n) \sim x \sum_{y < p \le z} \frac{1}{p-1}\,.
\end{equation}
By Markoff's inequality $H^*(x,y,z) \le M_1$, so it remains to prove the matching lower bound 
\begin{equation}\label{eq: thm short intervals PZ lower}
    H^*(x,y,z) \ge M_1 (1 + o(1))\,.
\end{equation}

First, suppose $z \le y + y^{3/5}$. 

Let us show that in this case \eqref{eq: thm short intervals PZ lower} holds regardless of whether we can compute $M_1$. Let $ M_2 = \sum_{n \leq x} r(n)^2$ and $M_2^\prime = M_2 - M_1 = \sum_{n \leq x} r(n)(r(n)-1)$, then by the Paley--Zygmund inequality 
\begin{equation}\label{eq:Markoff PZ}
    H^*(x,y,z) \ge \frac{M_1^2}{M_2} = \frac{M_1^2}{M_1 + M_2^\prime}\,.
\end{equation}
Bound the number of shifted prime pairs with the number of pairs, where only one of the divisors is required to be a shifted prime, then
\begin{equation}
    \begin{aligned}
        M_2^\prime = 2 \sum_{\substack{n \le x}} \sum_{\substack{y<p<q\le z \\ p-1,\,q-1|n}} 1 \le 2 \sum_{\substack{n \le x}} \sum_{\substack{y<p \le d\le z \\ p-1,\,d|n}} 1 \le
2x \sum_{\substack{y<p \le d\leq z}}\frac{(p-1,d)}{(p-1)d}\,.
    \end{aligned}
\end{equation}
Let $g = (p-1,d)$, then $g \le z-y \le y^{3/5}$, so 
\begin{equation}
    \begin{aligned}
       M_2^\prime  \le  2x \sum_{\substack{y<p \leq z }}\frac{1}{p-1} \sum_{\substack{g \mid (p-1)\\ g \le z-y}} \sum_{\substack{y < d \le z\\  g\mid d}} \frac{g}{d} & \ll x \sum_{\substack{y<p \leq z }}\frac{1}{p-1}   \sum_{\substack{g \mid (p-1)\\ g \le z-y}} \frac{z-y}{y} \,.
    \end{aligned}
\end{equation}
For any $n \le 2y$, we have $\tau(n) \ll y^{0.01}$, so
\begin{equation}
    \sum_{\substack{g \mid (p-1)\\ g \le z-y}} \frac{z-y}{y} \le y^{-2/5} \tau(p-1) \ll y^{-0.39}\,.
\end{equation}
Thus, we conclude that in this case $ M_2^\prime \ll y^{-0.39} M_1 = o(M_1)$ as required.

Now suppose $z \ge y + y^{3/5}$. 

Let $t = (z^\prime)^{1/2}$ and let
\begin{equation}\label{def: w=2loglogt}
    w:= 2\log\log t + \psi(t)\sqrt{\log\log t}, \quad \psi(t) = \min\{ (\log\log t)^{1/6}, \xi\}\,.
\end{equation}
Let
\begin{equation}
    \Tilde{M_1} := \sum_{n \leq x}r(n) \bm{1}_{\Omega(n,t) \le w} ,\quad \Tilde{M_2} := \sum_{n \leq x}  r(n)^2 \bm{1}_{\Omega(n,t) \le w}\,.
\end{equation}
By the Paley--Zygmund inequality, we have that
\begin{equation}\label{eq:restrictedlowerbound}
    H^*(x,y,z) \ge \sum_{n \le x} \bm{1}_{\exists p-1 \mid n, p \in (y,z]} \bm{1}_{\Omega(n,t) \le w} \ge \frac{\Tilde{M_1}^2}{\Tilde{M_2}} = \frac{\Tilde{M_1}}{1 + (\Tilde{M_2}-\Tilde{M_1})/\Tilde{M}_1}\,.
\end{equation}
By \Cref{prop:M1mild} we have $\tilde{M_1} = M_1 (1+ o(1))$. By \Cref{prop:M2 short intervals} for $y + y^{3/5} \le z \le 2y$
\begin{equation}\label{eq:thm short intervals M2 1}
        \tilde{M}_2^\prime = \tilde{M_2}-\tilde{M_1} \ll \tilde{M_1}\cdot (\log t)^{-(1/2 - \gamma_0)} \exp( 2 \psi(t) \sqrt{\log\log t}) = o(M_1)\,
\end{equation}
since $1/2-\gamma_0 = 3/2 - \log 4 \approx 0.11 > 0$. For $z \ge 2y$
\begin{equation}\label{eq:thm short intervals M2 2}
       \tilde{M}_2^\prime \ll \tilde{M_1}\cdot  (\log t)^{-(\gamma-\gamma_0)}  \exp( 2 \psi(t) \sqrt{\log\log t})\,.
\end{equation}
In this case the right-hand side is $o (M_1)$ as long as $\xi = (\gamma-\gamma_0) \sqrt{\log \log t} \ge 4 \psi(t) \to +\infty$. 

Thus, $ H^*(x,y,z) \ge \tilde{M}_1(1+ o(1)) = M_1(1+ o(1))$. This concludes the proof.
\end{proof}

Let us now present proofs of intermediate estimates required for \Cref{thm:short intervals}.

\subsection{First moment estimates}

We start with the following elementary lemma estimating the first moment of the representation function $r(n)$.

\begin{lemma}\label{lem:M1} Let $3 \leq y < z \leq x$ such that $z^\prime = \min\{x/y,z\} \to \infty$ as $x \to \infty$. Then as $x\to\infty$
\[
M_1 = \sum_{n \le x} r(n) = \left(x \sum_{p\in (y,z]} \frac{1}{p-1}\right) (1+ o(1))\,.
\]
In particular, for $y + y^{17/30+\varepsilon} \le z \le y^{1+o(1)}$, we have $M_1 \sim x u$, where $u=\frac{\log (z/y)}{\log z}$.
\end{lemma}
\begin{proof}
Write
\begin{equation}\label{eq proof manthm i M1 bound}
\begin{aligned}
        M_1 = \sum_{\substack{m(p-1) \leq x\\ p\in (y,z]}} 1 = 
        \sum_{p \in (y,z]} \sum_{m \leq \frac{x}{p-1}}1 
        =x \sum_{p \in (y,z]} \frac{1}{p-1} + O \brackets{ \sum_{p\in (y,z]} 1}\,.
\end{aligned}
\end{equation}
We claim that the first term in \eqref{eq proof manthm i M1 bound} always dominates.

Suppose first that $z \ge 2y$. In this case by Merten's estimate \cite[Theorem 23]{ingham1932distribution}
\begin{equation}
    \sum_{p \leq P}\frac{1}{p}=\log \log P + c_1 + O\left(e^{-c_2\sqrt{\log P}}\right)\,
\end{equation}
we have
\begin{equation}
   x \sum_{p \in (y,z]} \frac{1}{p-1} = x\log  \frac{\log z}{\log y } + O\left(xe^{-c_2\sqrt{\log y}}\right) =xu (1 + o(u)) + O\left(xe^{-c_2\sqrt{\log y}}\right)\,.
\end{equation}
By the prime number theorem the second term in \eqref{eq proof manthm i M1 bound} is of order
\begin{equation}
    \sum_{p\in (y,z]} 1 \asymp \frac{z}{\log z}  \asymp xu \cdot \frac{z/x}{\log (z/y)}
\end{equation}
If $\frac{x}{z}$ tends to $\infty$, then so does $\frac{x}{z} \log \frac{z}{y}$. If $x \asymp z$, then $x/y \ge z^\prime \to \infty$ implies that $z/y \to \infty$, so in this case $\frac{x}{z} \log \frac{z}{y} \to \infty$ as well. Thus, the Lemma holds for $z \ge 2y$.

Now assume $z \le 2y$. In this case, the first term in \eqref{eq proof manthm i M1 bound} is larger than the second by a factor of $x/z \ge x/(2y) \ge z^\prime/2 \to \infty$, so it indeed dominates. It remains to establish the asymptotics.

For $z \le 2y$ we are limited by existing estimates for counting primes in short intervals. In particular, for $y^{17/30+\varepsilon} \le z -y \le y$ by \cite{guth2026new} we have
\begin{equation}
    \sum_{p \in (y,z]} \frac{1}{p-1} = u\left(1+ O\left(u+ e^{-\sqrt[4]{\log y}}\right)\right)\,.
\end{equation}
 This concludes the proof of the Lemma.
\end{proof}

As is apparent from above, to compute $M_1$ one does not need to count prime factors of~$n$. However, the main contribution to $M_1$ indeed comes from integers with $\Omega^*(n,z^\prime) \approx 2 \log\log z^\prime$ prime factors up to $z^\prime = \min\left\{z,x/y\right\}$. Let us demonstrate this fact by showing that for any $t \le z^\prime$ the condition $\Omega(n,t) \le 2\log\log t + (\log\log t)^{1/2+\varepsilon}$ does not change the first moment.

\begin{lemma}\label{prop:M1mild}
    Let $3 \le y < z \le x$ such that $z \ge y + y^{17/30+\varepsilon}$ for some small fixed $\varepsilon>0$, and $z^\prime = \min\{z,x/y\} \to \infty$ as $x \to \infty$. Let $2 \le t \le z^\prime$. Let $w = 2\log\log t + \psi(t)\sqrt{\log\log t}$ for some positive function $\psi(t) \to +\infty$ such that $\psi = o(\sqrt{\log\log t})$. Then
    \[
   \sum_{n \leq x}r(n) \bm{1}_{\Omega(n,t) > w} \ll xu \cdot e^{-\frac{1}{6}\psi^2(t)}\,.    \]
\end{lemma}
\begin{proof}
To avoid dealing with sets of integers with fixed $\Omega(n,t)$ directly, we use Chernoff's inequality $\bm{1}_{\Omega\left(n,t\right)>w} \le \eta^{\Omega\left(n,t\right)-w}$ with
\begin{equation}
    \eta = \frac{w}{2\log\log t} = 1 + \frac{\psi(t)}{2\sqrt{\log\log t}} =: 1 + \epsilon_0 < 2\,
\end{equation}
so that
\begin{equation}
    \sum_{n\le x}{r(n)\bm{1}_{\Omega\left(n,t\right)>w}} = \sum_{y<p\le z}{1 \sum_{m\le\frac{x}{p-1}}{\bm{1}_{\Omega\left((p-1)m,t\right)>w}}} \le\eta^{-w}\sum_{y<p\le z}{\eta^{\Omega\left(p-1,t\right)}\sum_{m\le\frac{x}{p-1}}\eta^{\Omega\left(m,t\right)}}\,.
\end{equation}
By \Cref{lem:Tenenbaumsum} for the $m$-sum,
\begin{equation}
    \begin{aligned}
       \sum_{m\le\frac{x}{p-1}}\eta^{\Omega\left(m,t\right)} \ll  \frac{x}{p} (\log t)^{\eta -1} \,.
    \end{aligned}
\end{equation}

To bound the sum in $p$ we will use \Cref{lem:Koukoulopoulossum}. If $z \le 2y$, since we are in the short interval scenario, in order to apply \Cref{lem:Koukoulopoulossum} we need to ensure that $\epsilon_0 = o(1)$. This assumption is indeed satisfied given $\psi(t) = o(\sqrt{\log\log t})$, thus
\begin{equation}
    \sum_{y <p \leq z} \frac{\eta^{\Omega(p-1,t)}}{p} \ll \frac{z-y}{y}\cdot  \frac{ (\log \min\{t,z\})^{\eta-1}}{\log z} \ll u (\log t)^{\eta-1} \,.
\end{equation}

For $z \ge 2y$ we could have used \Cref{lem:Fordsieve} instead, which works for any $\eta < 2$ effectively increasing the range of admissible $\psi(t)$, but this will not be necessary for our purposes. Thus, in this case we have by \Cref{lem:Koukoulopoulossum}
\begin{equation}
    \begin{aligned}
        \sum_{y<p\le z}\frac{\eta^{\Omega\left(p-1,t\right)}}{p}
       &  \ll\sum_{\log{y}\le i\le\log{(4z)}}{\sum_{p\in I_i}\frac{\eta^{\Omega\left(p-1,t\right)}}{p}}
        \ll\left(\log{t}\right)^{\eta-1} \sum_{\log{y}\le i\le\log{(4z)}} \frac{1}{i} \ll u (\log t)^{\eta -1}\,
    \end{aligned}
\end{equation}
as well. We therefore deduce that
\begin{equation}
    \begin{aligned}
        \sum_{n\le x}{r(n)\bm{1}_{\Omega\left(n,t\right)>w}} \ll x \eta^{-w} \left(\log{t}\right)^{\eta-1}   \sum_{y<p\le z}\frac{\eta^{\Omega\left(p-1,t\right)}}{p}  \ll x u \cdot \eta^{-w} \left(\log{t}\right)^{2\left(\eta-1\right)}\,.
    \end{aligned}
\end{equation}
Using the Taylor expansion $\log(1 + \epsilon_0) = \epsilon_0 - \frac{\epsilon_0^2}{2}+O(\epsilon_0^3)$, we see that the exponent becomes 
\begin{equation}
    -w\log(1+\epsilon_0) + 2 \epsilon_0 \log \log t = -\epsilon_0^2 \log \log t +O(\epsilon_0^3 \log \log t) = \frac 14 \psi^2(t) (1+ o(1))\,.
\end{equation}
We obtain the strict upper bound $\ll xu  \cdot e^{-\psi^2(t)/6}$ as required.
\end{proof}

\subsection{Second moment estimate}

To prove the lower bound in \Cref{thm:short intervals} via the second moment method we require a bound on the average number of pairs of shifted prime divisors. The condition $\Omega(n,t) \le 2\log\log t + (\log\log t)^{1/2+\varepsilon}$ ensures that pairs of shifted prime divisors on $(y,z]$ almost never occur for $\gamma > \gamma_0$. Such a condition is unnecessary when $z \le y + o(y)$, but we include this range nonetheless for the sake of uniformity.

\begin{prop}\label{prop:M2 short intervals}
    Let $3 \le y < z \le x$ such that $z \ge y + y^{17/30+\varepsilon}$ for some small fixed $\varepsilon>0$, and $z^\prime = \min\{z,x/y\} \to \infty$ as $x \to \infty$. Let $(z^\prime)^{0.01} \le t\le z^\prime$. Let $w = 2\log\log t + g(t)$ for some positive function $g(t) \ge 5 \log\log\log t$. Then
\begin{equation}
 \Tilde{M}_2^\prime (w)  = \frac{1}{2} \sum_{n \leq x}  r(n)(r(n)-1) \bm{1}_{\Omega(n,t) \le w} \ll xu \cdot \exp( 2 g(t)) \cdot 
   \begin{cases}
        (\log t)^{-(1/2 - \gamma_0)},& \text{if $z\le 2y$}\,;\\
        (\log t)^{-(\gamma-\gamma_0)} ,& \text{if $z \ge 2y$}\,.
    \end{cases}
\end{equation}
\end{prop}
\begin{proof}
To avoid working with prime factors of $n$ directly, we use Chernoff's inequality with exponent $\eta = 2$ so that $\bm{1}_{\Omega\left(n,t\right)\le w} \le 2^{w-\Omega(n,t)}$ and therefore
    \begin{equation}
    \begin{aligned}
\Tilde{M}_2^\prime
      = \sum_{\substack{n \le x}} \bm{1}_{\Omega\left(n,t\right)\le w}\sum_{\substack{y<p<q\le z \\ p-1,\,q-1|n}} 1 & =   \sum_{\substack{y<p<q\le z\\
        \ell = [p-1,q-1]\le x}}  \sum_{m\le\frac{x}{\ell}}1 \\
        & \le 2^{w}  \sum_{\substack{y<p<q\le z\\
        \ell = [p-1,q-1]\le x}} 2^{-\Omega\left(\ell,t\right)}  \sum_{m\le\frac{x}{\ell}}2^{-\Omega\left(m,t\right)} \,. 
    \end{aligned}
    \end{equation}
Set $L(h) := L(h;x,t) =  \log(2 + \min\{x/h,t\})$. Then by \Cref{lem:Tenenbaumsum} for the $m$-sum, where $m \le x/\ell \le x/y$, we have that
\begin{equation}\label{eq:propM2 regime 1 S}
    \begin{aligned}
        \Tilde{M}_2^\prime
        \ll x2^{w}  \sum_{\substack{y<p<q\le z\\
        \ell = [p-1,q-1]\le x}} \frac{2^{-\Omega\left(\ell,t\right)}}{\ell }\cdot L(\ell)^{-1/2}\,.
    \end{aligned}
\end{equation}

Let $g=(p-1,q-1)$, then $p-1 = ga$ and $q-1=gb$ for some coprime integers $a < b$, and $\ell = gab$. Since $y < ga < gb \le z$ and $g\le q-p \le z-y$, we have 
\begin{equation}\label{eq:ab range}
 \frac{y}{z-y} \le a < b\le \min\{z/g,x/p\} \le \min\{z,x/y\} = z^\prime\,.
\end{equation}
On the other hand, $\ell \le x$ and $\ell = (p-1)(q-1)/g \ge y^2/g$, so $g \ge G_0 := \max \{y^2/x, 1\}$. Hence $g$ must lie in the range $G_0 \le g \le z-y$.

To count $a,b$ and $g$ with two primality conditions we use sieve methods. When $g$ is small, say $g \le G$, we fix $g$ and use relations $p,q\equiv1(g)$ to count $a$ and $b$ separately, one primality condition per variable. When $g$ is large, say $g \ge G$, we instead fix $a$ and $b$ and count $g$ in arithmetic progressions $g\equiv1(a)$ and $g\equiv1(b)$ satisfying two primality conditions simultaneously. The latter range is more challenging and always has a bigger contribution. 

The cut-off $g\le G$ is relatively flexible. We can assume that $G_0  y^{0.1}~\le~G~\le~\min\{z-y, y^{0.6}\}$, which is only meaningful if $G_0 \le y^{0.5}$ and therefore $y \le x^{2/3}$. With this in mind, we set $G = G_0 y^{0.1}$ and split $\Tilde{M}_2^\prime$ into two parts such that
\begin{equation}\label{eq:prop3 S < S1 + S2}
    \Tilde{M}_2^\prime \ll x 2^w (S_1 + S_2)\,,
\end{equation}
where $S_1$ is the part of the sum in \eqref{eq:propM2 regime 1 S} with $G_0 \le g \le G$, and $S_2$ is the part with $g \ge G$. 

Let us first show that
\begin{equation} \label{eq: prop3 S1 final}
    S_1 \ll \frac{u^2}{\log t}\log\log G\,.
\end{equation}
For this part we can assume that $y \le x^{2/3}$, which also implies that $\log z^\prime \asymp \log z \asymp \log G$ and $\log z^\prime \ge \log G_0$. Write $\ell = [p-1,q-1] = (p-1)(q-1)/g$, then
\begin{equation}\label{eq:S1 def}
\begin{aligned}
S_1 & \le \sum_{G_0\le g\le G} g2^{-\Omega\left(g,t\right)} \sum_{\substack{y<p,q\le z \\ p,q\equiv1(g)}} \frac{2^{-\Omega\left(\frac{p-1}{g},t\right)-\Omega\left(\frac{q-1}{g},t\right)}}{pq \sqrt{L\left(pq/g\right)}}  =: \sum_{G_0\le g\le G} \frac{2^{-\Omega\left(g,t\right)}}{g} s_1(g)\,.
\end{aligned}
\end{equation}

If $z \le 2y$, we have $xg/pq \ge xg/z^2 \ge xg/(4y^2)$ and thus $L\left(pq/g\right) \asymp L\left(y^2/g\right)$. By \Cref{lem:Koukoulopoulossum} for $h\le P/2$ together with \Cref{lem:Tenenbaumsum} for $f(n) = 2^{-\Omega(n,t)}$ we have
    \begin{equation}\label{eq: prop2 sum over p}
        \sum_{\substack{P-h<p\le P \\ p\equiv1(g)}} \frac{2^{-\Omega\left(\frac{p-1}{g},t\right)}}{p} \ll \frac{1}{P}\cdot \frac{h}{\varphi(g)\log P} \cdot (\log t)^{-1/2}\,.
    \end{equation}
Taking $P = z$ and $h = z-y$, such that the above turns into $\ll \frac{u}{\varphi(g) \sqrt{\log t}}$, we conclude that
    \begin{equation}\label{eq:prop2 s_1(g) zge2y}
        \begin{aligned}
            s_1(g) \ll g^2 \brackets{\frac{u}{\varphi(g)\sqrt{\log t} }}^2 L\left(y^2/g\right)^{-1/2} \ll
            \frac{u^2}{\log t} \left(\frac{g}{\phi(g)} \right)^2  L\left(y^2/g\right)^{-1/2}\,.
        \end{aligned}
    \end{equation}
If $z\ge 2y$, fix the size of $p$ and $q$ and split the sum into intervals $p \in I_i = (e^{i-1},e^{i}]$ and $q \in I_j = (e^{j-1},e^{j}]$. Using \eqref{eq: prop2 sum over p}, we reduce $s_1(g)$ to     
    \begin{equation}\label{eq:s_1g zge2y 1}
        \begin{aligned}
            s_1(g) & \ll g^2\sum_{\log y < i,j\le \log (4z)} L\left(e^{i+j}/g\right)^{-1/2} \sum_{\substack{p \in I_i,\, q \in I_j \\ p,q\equiv1(g)}} \frac{2^{-\Omega\left(\frac{p-1}{g},t\right)-\Omega\left(\frac{q-1}{g},t\right)}}{pq} \\
            & \ll \frac{g^2(\log t)^{-1}}{\varphi^2(g)(\log z)^2}\sum_{\log y < i,j\le \log (4z)} L\left(e^{i+j}/g\right)^{-1/2}\,.
        \end{aligned}
    \end{equation}
Observe that $\sum_{A < k \le B} k^{-c} \ll_c (B-A) B^{-c}$ for any $c \in [0,1)$ and positive integers $A < B$, whence for any $x_1,t_1 > 0$
\begin{equation}
\label{eq:s_1g zge2y 2}
    \begin{aligned}
                % \sum_{A < j \le B} L(e^j/C)^{-1/2} & = 
                \sum_{A < j \le B} \left(\log (2 + \min\{x_1 e^{-j},t_1\})\right)^{-1/2} \ll (B-A) \cdot \left(\log (2 + \min\{x_1 e^{-A},t_1\})\right)^{-1/2}\,.
    \end{aligned}
\end{equation}
Applying \eqref{eq:s_1g zge2y 2} to the final sum in \eqref{eq:s_1g zge2y 1} we see that 
\begin{equation}
    \begin{aligned}
        \sum_{\log y < i,j\le \log (4z)} L\left(e^{i+j}/g\right)^{-1/2} & \ll \log (z/y) \sum_{2\log y < j\le 2\log (4z)} L\left(e^{j}/g\right)^{-1/2} \\
        & \ll  \log^2 (z/y) L\left(y^2/g\right)^{-1/2}\,,
    \end{aligned}
\end{equation}
and thus verify that \eqref{eq:prop2 s_1(g) zge2y} also holds for $z \ge 2y$.

To bound $S_1$ it now remains to evaluate the sum over \eqref{eq:prop2 s_1(g) zge2y}:
\begin{equation}\label{eq: prop3 S1 final 0}
    \begin{aligned}
         S_1 \ll \frac{u^2}{\log t} \sum_{\log G_0 < i\le \log (4G)} L\left(y^2e^{-i}\right)^{-1/2} \sum_{g \in I_i}  \frac{1}{g} \left(\frac{g}{\varphi(g)}\right)^2 2^{-\Omega\left(g,t\right)}\,.
    \end{aligned}
\end{equation}
By \Cref{lem:Tenenbaumsum} for $g \in I_i = (e^{i-1},e^{i}]$, and any $k\in \N$, we have
\begin{equation}\label{eq: sum 2^-Omega}
    \sum_{g \in I_i}  \frac{1}{g} \left(\frac{g}{\phi(g)}\right)^k 2^{-\Omega\left(g,t\right)} \ll_k i^{-1/2}\,,
\end{equation}
so
\begin{equation}
    \begin{aligned}
         L\left(y^2e^{-i}\right)^{-1/2} \sum_{g \in I_i}  \frac{1}{g} \left(\frac{g}{\varphi(g)}\right)^2 2^{-\Omega\left(g,t\right)} \ll L\left(y^2e^{-i}\right)^{-1/2}i^{-1/2} =: F_1(i)\,.
    \end{aligned}
\end{equation}
For any $A, B \ge 1$ we have $(AB)^{-1} \ll A^{-2} + B^{-2}$, and $\min\{A,B\}^{-1} \ll A^{-1} + B^{-1}$. Then for $i \le \log (4G) \ll \log z$ the above satisfies
\begin{equation}
\label{eq:prop3 F1(i) bound}
\begin{aligned}
        F_1(i) 
        % \ll i^{-1} + L\left(y^2e^{-i}\right)^{-1}
         \ll i^{-1} + (\log (2 + xe^i/y^2))^{-1}\,.
\end{aligned}
\end{equation}
Summing over $\log G_0 < i\le \log (4G)$ yields 
\begin{equation} \label{def:prop3 S2}
    \begin{aligned}
        \sum_{\log G_0 < i\le \log 4G} F_1(i) \ll \sum_{\log G_0 <i\le \log 4G}\brackets{ i^{-1} + \left(\log{\left(2+xe^i/y^2\right)}\right)^{-1}} \ll \log\log G
    \end{aligned}
\end{equation}
and therefore \eqref{eq: prop3 S1 final}.

Let us now turn to $S_2$. This part is always present in the sum, but we cannot assume that $\log z^\prime \asymp \log z$ anymore. We will show that for $z\le 2y$
\begin{equation}\label{eq: prop3 S2 final zle2y}
    \begin{aligned}
        S_2 \ll \frac{u}{(\log z) \sqrt{\log t}} \cdot (\log\log t)^5 \ll 2^{-w} \cdot u (\log t)^{-(\frac{1}{2} - \gamma_0)} \exp(2g(t))\,,
    \end{aligned}
\end{equation}
and for $z \ge 2y$
\begin{equation}\label{eq: prop3 S2 final zge2y}
   S_2 \ll \frac{u^2}{\log{t}} \cdot (\log\log t)^5 \ll 2^{-w} \cdot u (\log t)^{-(\gamma-\gamma_0)}  \exp( 2 g(t))\,.
\end{equation}
Both bounds are larger than the bound for $S_1$ in \eqref{eq: prop3 S1 final}. The bound for $z \le 2y$ is weaker, however, due to the length of the interval, we can afford to be relatively crude, so even though \eqref{eq: prop3 S2 final zle2y} is not optimal, it will suffice.

Observe that $y < ga < gb \le z $ implies that $a < b < a z/y$. Then setting $t_2 = \min\{z^\prime,z/G\}$ and using \eqref{eq:ab range} we can bound $S_2$ by
\begin{equation}
    \begin{aligned}
        S_2  & \le \sum_{\substack{\frac{y}{z-y} < a < b \le t_2\\
         a < b < a \frac{z}{y}}} \frac{2^{-\Omega(ab,t)}}{ab} \sum_{\substack{\frac{y}{b} < g\le \frac{z}{b}\\
        ag+1,bg+1 \in \PP}} \frac{2^{-\Omega(g,t)}}{g\sqrt{L(abg)}}  =: \sum_{\substack{\frac{y}{z-y} < a < b \le t_2\\
         a < b < a \frac{z}{y}}} \frac{2^{-\Omega(ab,t)}}{ab} s_2(a,b)\,.
    \end{aligned}
\end{equation}

To bound $s_2(a,b)$, split the sum into intervals $g \in I_i = (e^{i-1},e^{i}]$; if $z \le 2y$, complete the interval such that
\begin{equation}
    s_2(a,b) 
         \ll \sum_{\log{\frac{y}{b} < i\le\log{\frac{4z}{b}}}}  L(abe^i)^{-1/2} \sum_{\substack{g \in I_i\\
        ag+1,bg+1 \in \PP}} \frac{2^{-\Omega(g,t)}}{g}\,.
\end{equation}
By Lemma \ref{lem:Fordsieve} together with the inequality $N/\varphi(N) \ll \log\log N$, for any $i \gg \log z$ we have
\begin{equation}\label{eq:prop3 S2 zle 2y 1}
\begin{aligned}
        \sum_{\substack{g \in I_i\\
        ag+1,bg+1 \in \PP}} \frac{2^{-\Omega(g,t)}}{g} \ll \frac{\Phi(a,b)}{(\log z)^2 \sqrt{\log t}}  \ll \frac{(\log\log t)^4 }{(\log z)^2 \sqrt{\log t}}\,.
\end{aligned}
\end{equation}    
Apply \eqref{eq:prop3 S2 zle 2y 1} together with \eqref{eq:s_1g zge2y 2}, then uniformly in $b$
    \begin{equation}\label{prop2-16}
        \begin{aligned}
        s_2(a,b) 
        & \ll \frac{(\log\log t)^4 }{(\log z)^2 \sqrt{\log t}} \sum_{\log{\frac{y}{b}< i\le\log{\frac{4z}{b}}}} L(abe^i)^{-1/2} \\
        & \ll \frac{(\log\log t)^4 \log (4z/y)}{(\log z)^2 \sqrt{\log t}}  \cdot L(ya)^{-1/2}\,.
        \end{aligned}
    \end{equation}

Suppose that $z\le 2y$. 

In this case we can crudely bound the sum in $b$ by
\begin{equation}\label{eq:prop3 S2 zle 2y 2}
    \begin{aligned}
        \sum_{a < b < az/y} \frac{2^{-\Omega(b,t)}}{b} \le \sum_{a < b < az/y} \frac{1}{b} \le \frac{1}{a} \sum_{a < b < a + a(z-y)/y} 1 \le \frac{z-y}{y} \ll \log (z/y)\,.
    \end{aligned}
\end{equation}
Insert \eqref{prop2-16} and \eqref{eq:prop3 S2 zle 2y 2} into \eqref{def:prop3 S2}, then indeed as claimed in \eqref{eq: prop3 S2 final zle2y}
\begin{equation}\label{eq:prop3 S2 short 1}
    \begin{aligned}
        S_2 & \ll  \frac{(\log\log t)^4\log (z/y)}{(\log z)^2\sqrt{\log t}} \sum_{a\le t} \frac{2^{-\Omega(a)}}{a} L(ya)^{-1/2} \ll \frac{(\log\log t)^4\log (z/y)}{(\log z)^2\sqrt{\log t}} \cdot \log\log t\,.\,
    \end{aligned}
\end{equation}
The last step (evaluating the sum in $a$) is done similarly to \eqref{eq: prop3 S1 final 0} by applying \eqref{eq: sum 2^-Omega} to each interval $a \in I_i = (e^{i-1},e^{i}]$ and and invoking \eqref{eq:prop3 F1(i) bound} for the final sum in $i\le \log (4t)$.

Now suppose $z \ge 2y$. 

To bound the sum in $b$, split into $e$-adic intervals and apply \eqref{eq: sum 2^-Omega}, then
\begin{equation}\label{eq:prop 3 S2 bsum}
    \begin{aligned}
        \sum_{a < b < \frac{az}{y}} \frac{2^{-\Omega(b)}}{b} & \ll \sum_{\log a < i \le \log \frac{4az}{y}} \sum_{b \in I_i} \frac{2^{-\Omega(b)}}{b}  \ll \sum_{\log a < i \le \log \frac{4az}{y}} i^{-1/2} \ll \frac{\log (z/y)}{(\log (1 + a))^{1/2}}\,.
    \end{aligned}
\end{equation}
Combining \eqref{prop2-16} and \eqref{eq:prop 3 S2 bsum} we see that it now remains to evaluate the sum over $a$:
\begin{equation}\label{eq:prop3 S2 1}
    \begin{aligned}
         S_2 & \ll \frac{u^2(\log\log t)^4}{\sqrt{\log t}} \sum_{\substack{a \le t}} \frac{2^{-\Omega(a)}}{a} L(ya)^{-1/2} (\log (1 + a))^{-1/2}\,.
    \end{aligned}
\end{equation}
Comparing the sum above to \eqref{eq: prop3 S1 final 0} we observe serve that 
\begin{equation}\label{eq:prop3 S2 2}
    \begin{aligned}
        \sum_{\substack{a \le t}} \frac{2^{-\Omega(a)}}{a} L(ya)^{-1/2} (\log (1 + a))^{-1/2} & \ll \sum_{\substack{1 \le i \le \log (4t)}} F_1(i) i^{-1/2}  =: \sum_{\substack{1 \le i \le \log (4t)}} F_2(i)\,,
    \end{aligned}
\end{equation}
where $F_2(i) := L(ye^i)^{-1/2} i^{-1}$.

If $i \le \frac12 \log t$, then $L(ye^i) \gg L((xy)^{1/2}) \asymp \log t$, so $F_2(i) \ll (\log t)^{-1/2} i^{-1}$. If $i \ge \frac12 \log t$, then $F_2(i) \ll  L(ye^i)^{-1/2} (\log t)^{-1}$. Applying \eqref{eq:s_1g zge2y 2} in the latter case, we obtain 
\begin{equation}\label{eq:prop3 S2 3}
    \begin{aligned}
        \sum_{\substack{1 \le i \le \log (4t)}} F_2(i) & \ll  \frac{1}{\sqrt{\log t}} \sum_{\substack{1 \le i \le \frac12 \log t}} \frac{1}{i} + \frac{1}{\log t} \sum_{\substack{\frac12 \log t \le i \le \log (4t)}} L(ye^i)^{-1/2}\\
        & \ll \frac{\log\log t}{\sqrt{\log t}} + L((xy)^{1/2})^{-1/2} \ll \frac{\log\log t}{\sqrt{\log t}}\,.
    \end{aligned}
\end{equation}

Combining the above with \eqref{eq:prop3 S2 1} and \eqref{eq:prop3 S2 3} yields the desired bound \eqref{eq: prop3 S2 final zge2y}.

 This concludes the proof of the Proposition.
\end{proof}

%% file: 05_long_intervals.tex
In this section we prove parts (ii), (iii) and (iv) of \Cref{mainthm} for the remaining regime where $0 < \gamma\le \gamma_0 + 0.01$. For $\gamma > 0$ we establish the order of magnitude of $H^*(x,y,z)$. When $\gamma = o(1)$, our lower bound fails to match the upper bound, so in this case we only establish the order of magnitude of $H^*(x,y,z)$ up to a factor of $(\log z^\prime)^{o(1)}$. 

\begin{proof}[Proof of \Cref{thm:longer intervals} and \ref{thm:gamma to 0}]
Let $z \ge 2y$. Choose any $t$ in the range $(z^\prime)^c \le t \le z^\prime$. Let $\lambda = \log\log t$. Set $\rho:=\frac{1+\gamma}{\log 2}$, $k_0 := \lceil \rho\log\log{t}\rceil$ and $w:=2\log\log t + 6(\log\log t)^{2/3}$. 

Let $H^*_k (x,y,z) = \#\{n \le x: \Omega(n,t) = k,\ \exists p \in (y,z],\ p-1 \mid n\}$, then
\begin{equation}\label{thm:long intervals H = sum Hk}
    \begin{aligned}
         H^*(x,y,z) & = \sum_{k \ge 0} H^*_k(x,y,z)  = \sum_{k \le w} H^*_k(x,y,z) + \sum_{\substack{n\le x\\ \Omega(n,t) > w}} r(n)\,.
    \end{aligned}
\end{equation}
By \Cref{prop:M1mild} 
\begin{equation}\label{eq:thm long intervals k>w}
    \sum_{\substack{n\le x}} r(n) \bm{1}_{\Omega(n,t) > w} \ll xu \exp(-(\log\log t)^{1/3})\,,
\end{equation}
so bounding \eqref{thm:long intervals H = sum Hk} boils down to bounding $H^*_k(x,y,z)$ for each $k\le w$.

Let us first bound $H^*(x,y,z)$ from above.

Let $\pi_k = \#\left\{n \le x : \Omega(n,t)=k\right\}$, and $M_{1,k}=\sum_{\substack{n \leq x}}r(n)\bm{1}_{\Omega(n,t)=k}$. By Lemma \ref{lem:Big-Omega-bounds}
\begin{equation}
    \pi_k \ll \pi_k^\prime := \frac{x(\log\log t)^k}{k!\log t} = \frac{x\lambda^k}{k!\log t}\,,
\end{equation}
and by Lemma~\ref{lem: M1 fixed k} for $k\le w$ we have 
\begin{equation}
    M_{1,k} \asymp \frac{2^ku}{\log t} \cdot \pi_k^\prime =  \frac{2^k}{(\log t)^{1+\gamma}}\cdot \pi_k^\prime\,.
\end{equation}
If $k_0 \ge 2\log\log t$ and therefore $\gamma \ge \gamma_0$, it suffices to use Markoff's inequality so that
\begin{equation}\label{eq:thm long intervals rho>2}
    \begin{aligned}
        \sum_{k \le w} H^*_k(x,y,z) & \le \sum_{k\le w} M_{1,k} \asymp \frac{x}{\left(\log{t}\right)^{\gamma}}\sum_{k\le w}\frac{\left(2\lambda\right)^k}{k!e^{2\lambda}}  \asymp \frac{x}{(\log t)^\gamma} = xu\,.
    \end{aligned}
\end{equation}
If $k_0 < 2\log\log t$, then for each $k \le w$ take the minimum between Markoff's inequality and the trivial bound:
\begin{equation}
    H^*_k(x,y,z) \le  \min \{M_{1,k}, \pi_k\} \ll  \min\left\{M_{1,k}, \pi_k^\prime\right\}\,.
\end{equation}
Observe that $M_{1,k} \ll \pi_k^\prime$ for $k \le k_0$, and $M_{1,k} \gg \pi_k^\prime$ for $k > k_0$. Then, applying \Cref{lem:improvedBenett}, we see that
\begin{equation}\label{eq:thm long intervals rho<2}
    \begin{aligned}
    \sum_{k \le w} H^*_k(x,y,z) & \le \sum_{k\le k_0} M_{1,k}+\sum_{k_0<k\le w}\pi_k \ll \frac{x}{\left(\log{t}\right)^{\gamma}}\sum_{k\le k_0}\frac{\left(2\lambda\right)^k}{k!e^{2\lambda}} + x\sum_{k>k_0}\frac{\lambda^k}{k!e^\lambda}\\
    & \ll \frac{x}{(\log t)^{Q(\rho)}((\gamma_0-\gamma) \lambda^{1/2} + 1)} + \frac{x}{(\log t)^{Q(\rho)}\lambda^{1/2}}\,,
    \end{aligned}
\end{equation} 
where the first term always dominates the bound. Gathering \eqref{eq:thm long intervals k>w}, \eqref{eq:thm long intervals rho>2} and \eqref{eq:thm long intervals rho<2} we obtain the desired upper bound for $H^*(x,y,z)$.

Let us now move on to the lower bound. If $\gamma\ge 0.02$, use the sum
\begin{equation}
    H^*(x,y,z) \ge \sum_{k_0-0.01\lambda \le k \le k_0} H^*_{k}(x,y,z)\,.
\end{equation}
If $\gamma < 0.02$, it suffices to use the lower bound
\begin{equation}
     H^*(x,y,z) \ge H^*_{k_0}(x,y,z)\,.
\end{equation}
In either case, we can assume that $k \ge \frac{1}{\log 2} \log \log t$.

First, we would like to restrict divisors of $n$ to a subset of positive proportion to avoid certain technicalities while working with shifted primes. If $n = (p-1)d\le x$, we assume that $x/(10z) \le d \le x/y$, and that both $p-1$ and $d$ have their respective largest prime factors satisfy $P^+(p-1) \ge (p-1)^c$, $P^+(d) \ge d^c$ for some small fixed $c > 0$. That is, define
\begin{equation}
    \tilde{r} (n) = \#\{(p-1)d= n: p\in (y,z],\ d\in (x/z,x/y],\ P^+(p-1)\ge z^c, P^+(d) \ge (x/y)^c\}\,.
\end{equation}
Next, we would like to restrict $n$ to a subset of positive proportion within the set of integers with $\Omega(n,t) = k$. We assume that $P^+(n) \ge x^c$, and that $n$ satisfies the following barrier condition.

Let $J = \lceil\log\log t\rceil$ and $t_j = e^{e^{j}}$ for $1\le j \le J$. Let $J_0 = \lfloor\log \log (z/y)\rfloor$.

Let $\{B(j)\}_{j=1}^J$ be a set of positive constants such that 
\begin{equation}\label{eq: thm long intervals barrier local}
    B(j) = \rho j + \rho \cdot g(\min\{j,J-j\}),\quad g(s) = s^{1/2 + \varepsilon}\,.
\end{equation}
Let $C \ge 100$ be large and fixed, and let
    \[
    \mathcal{B}(k) = \{n\le x:  P^+(n) \ge x^c, \ \Omega(n,t) = k,\ \Omega(n,t_j) \le B(j)+C\ \forall 1\le j \le J\}\,.
    \]
Then 
\begin{equation}
    H^*_k(x,y,z) \ge \sum_{\substack{n\le x\\ \Omega(n,t)=k}} \bm{1}_{\tilde{r} (n) \ge 1} \ge \sum_{\substack{n\le x}} \bm{1}_{\tilde{r} (n) \ge 1} \bm{1}_{n \in \mathcal{B}(k)} =: \tilde{H}^*_k(x,y,z) \,.
\end{equation}
Let $\tilde{M}_{1,k}(\mathcal{B})= \sum_{\substack{n\leq x}}\tilde{r}(n)\bm{1}_{n \in \mathcal{B}(k)}$, $\tilde{M}_{2,k}(\mathcal{B})= \sum_{\substack{n\leq x}}\tilde{r}^2(n)\bm{1}_{n \in \mathcal{B}(k)}$, then by the Paley--Zygmund inequality
\begin{equation}\label{eg: thm long intervals lower bound 1}
    \tilde{H}^*_k(x,y,z) \ge \frac{\tilde{M}^2_{1,k}(\mathcal{B})}{\tilde{M}_{2,k}(\mathcal{B})} \,.
\end{equation}

By \Cref{prop:firstmomentbarrier} we have $\tilde{M}_{1,k}(\mathcal{B}) \asymp M_{1,k}$ for sufficiently small $c>0$. By \Cref{prop:second moment general barrier}
\begin{equation}
    \tilde{M}_{2,k}^\prime(\mathcal{B}) = \tilde{M}_{2,k}(\mathcal{B}) - \tilde{M}_{1,k}(\mathcal{B}) \ll M_{1,k} V_k\,,
\end{equation}
where
\begin{equation}
    V_k =  u \sum_{j\le J} 2^{B(j)} (e^{j} + \log (z/y))^{-1} \ll u \sum_{J_0 \le j\le J} 2^{B(j)}e^{-j}\,
\end{equation}
since
\begin{equation}
    \begin{aligned}
        \sum_{j\le J_0} 2^{B(j)} & \ll  \sum_{j\le J_0} 2^{\rho (j + g(\min\{j,J-j\}))} \ll 2^{\rho (J_0 + g(\min\{J_0,J-J_0\}))} 
    \end{aligned}
\end{equation}
is dominated by the terms at $j = J_0 +O(1)$.

Let $f(j) = -(\rho\log 2 - 1)j + \rho g(\min\{j,J-j\})$, then 
\begin{equation}
    B(j) \log 2 - j = (\rho\log 2 - 1)j + \rho g(\min\{j,J-j\}) = J(\rho\log 2 - 1) -f(J-j)
\end{equation}
so that
 \begin{equation}
 \begin{aligned}
          V_k &  \ll u\sum_{J_0 \le j\le J} e^{(\rho\log 2 - 1)j + \rho g(\min\{j,J-j\})} \ll  \frac{2^ku}{\log t} \sum_{j\le J-J_0} e^{-f(j)}\,.
 \end{aligned}
 \end{equation}
If $\gamma \ge 0.02$, then $k \ge k_0-0.01\lambda$ implies that $\rho\log 2 - 1 \ge \gamma -0.01\ge 0.01$, so 
\begin{equation}
    \sum_{j\le J-J_0} e^{-f(j)} \le \sum_{j\le J-J_0} e^{-0.01j + \rho g(\min\{j,J-j\})} \ll 1 \,.
\end{equation}
Thus, for all such $k$ we have 
\begin{equation}
    \frac{ \tilde{M}_{2,k}(\mathcal{B})}{ \tilde{M}_{1,k}(\mathcal{B})} \ll 1 + \frac{2^k u}{\log t}\ll 1 + \frac{2^{k_0} u}{\log t} \ll 1\,.
\end{equation}
Summing \eqref{eg: thm long intervals lower bound 1} over $k_0 - 0.01 \le k\le k_0$ yields the desired lower bound matching the upper bound \eqref{eq:thm long intervals rho<2}. 

Now suppose $\gamma \le 0.02$, so we can take $k = k_0$. In this case without the loss of generality we can assume that $J_0 \ge J/2$ so that $f(j) = -(\rho\log 2 - 1)j + \rho g(j) = -(\rho\log 2 - 1)j + \rho j^{1/2+\varepsilon}$. 

Let $j_0 = \lceil(2\rho/\gamma)^{\frac{1}{1/2-\varepsilon}}\rceil$ be the largest integer such that $\gamma j \le 2\rho g(j)$, then 
\begin{equation}
    \begin{aligned}
        V_k  \le \sum_{\substack{j\le J-J_0\\ j \le j_0}} e^{\rho j^{1/2+\varepsilon}} + \sum_{\substack{j\le J-J_0\\ j\ge j_0 }} e^{-\gamma j/2} & \ll \exp \brackets{C_1 \min\{J-J_0, j_0\}^{1/2+\varepsilon}}\\
        & \ll \exp \brackets{C_2\min\{J-J_0, 1/\gamma^2\}^{1/2+2\varepsilon}}
    \end{aligned}
\end{equation}
for some positive constants $C_1, C_2 > 0$ depending on $\varepsilon$. If $\gamma > 0$ or $J-J_0 \ll 1$, this factor is of size $1$ as required. This concludes the proof of \Cref{thm:longer intervals} and \Cref{thm:gamma to 0}.
\end{proof}

We now prove the intermediate statements required to complete the above proof.

\subsection{First moment estimates}

We start by proving a technical lemma for the first moment of the representation function $r(n)$ twisted by a multiplicative function.

\begin{lemma} \label{lem:first moment multiplicative}
Let $3 \le y < z\le x$ such that $z \ge 2y$ and $z^\prime = \min\{z, x/y\} \to \infty$ as $x \to \infty$. Let $f \in \mathcal{F}_{b,\eta}$ for some $b \in \N$ and $\eta \in [1,2)$. Let $z/y \le t\le z^\prime$. Let $\lambda(f) = \sum_{p \leq t}f(p)/p$. For $\ell \in \N$ define
\begin{equation}
    \pi_\ell(x;f) = \frac{x}{\ell!\log x} \brackets{\sum_{p \le t} \frac{f(p)}{p}}^\ell \exp \brackets{\sum_{t<p\le x} \frac{f(p)}{p}}\,.
\end{equation}
% Let $t^\prime$ be an arbitrary number such that $3 \le t^\prime \le t$.
Then uniformly for integers $k$ such that $0\le k \le \frac{4-\varepsilon}{\eta} \lambda(f)$ we have  
        \begin{equation}
        M_{1,k}(f) = \sum_{\substack{n\leq x\\
    \Omega(n,t) = k}}r(n) f(n) \ll_{b,\eta,\varepsilon} \frac{2^ku}{\log z^\prime} \exp \brackets{\sum_{t<p\le z^\prime} \frac{f(p)}{p}}\cdot \pi_{k}(x;f)\,.
    \end{equation}
\end{lemma}
\begin{proof}
Let $\rho = k/\lambda(f) \le \frac{4-\varepsilon}{\eta}$. Let $k_0 = \frac{2-\varepsilon/2}{\eta} \lambda(f)$. 

Let $n = (p-1)m \le x$ such that $\Omega(n,t) = \Omega(p-1,t) + \Omega(m,t) = k$. Fix $k_1 = \Omega(p-1,t)$ and $k_2 = \Omega(m,t)$. The main contribution to $M_{1,k}(f)$ comes from
\begin{equation}
    k_1 \approx k_2 \approx k/2 \le  \frac{2-\varepsilon/2}{\eta} \lambda(f)\,,
\end{equation}
which is well within the admissible range for Poisson estimates. The contribution from $k_i \ge (\rho/2+\varepsilon_0)\lambda(f)$ is small in comparison and does not require fixing $k_i$. To see this, let us first estimate
\begin{equation}\label{lem: M1 fixed k-1}
        S(k_1,k_2) = \sum_{\substack{y <p \le z \\ \Omega \left(p-1,t\right)=k_1}} \sum_{\substack{m\le x/p \\ \Omega\left(m,t\right)=k_2}} f(m(p-1)) = \sum_{\substack{m\le x/y \\ \Omega\left(m,t\right)=k_2}}\sum_{\substack{y<p \le \min\{z,x/m\} \\ \Omega \left(p-1,t\right)=k_1}}  f(m(p-1))
    \end{equation}  
for $k_i \le k_0$ via \Cref{lem: hardy-ramanujan multiplicative} and \Cref{lem: hardy-ramanujan multiplicative shifted prime}.

Suppose $\log (x/y) \le 100 \log( z/y) \le 100\log z^\prime$, then $\log (z/y) \asymp \log z^\prime \asymp \log (x/y)$ and $\log z \asymp \log x$. Sum over $p$ first, then
\begin{equation}
    \begin{aligned}
        S(k_1,k_2) & \le \sum_{\substack{m\le x/y \\ \Omega\left(m,t\right)=k_2}} f(m) \sum_{\substack{y<p \le x/m \\ \Omega \left(p-1,t\right)=k_1}}  f(p-1) 
        \\
        & \ll_{b,\eta}  \frac{x}{(\log z)^2} \exp \brackets{\sum_{t<p\le z} \frac{f(p)}{p}} \frac{\lambda(f)^{k_1}}{k_1!} \sum_{\substack{m\le x/y \\ \Omega\left(m,t\right)=k_2}}  \frac{f(m)}{m}  \\
       &  \ll_{b,\eta}  \frac{x}{(\log z)^2}\exp \brackets{\sum_{t<p\le z} \frac{f(p)}{p} + \sum_{t<p\le x/y} \frac{f(p)}{p}}
       \cdot \frac{\lambda(f)^{k_1}\lambda(f)^{k_2}}{k_1!k_2!}\,.
    \end{aligned}
\end{equation}
If $\log (x/y) \ge 100 \log (z/y)$, then $\log (x/z) = \log (x/y) -\log (z/y) \ge 0.99 \log (x/y)$ and $ \log (x/z) \le \log (x/p) \le \log (x/y)$, so $\log (x/p) \asymp \log (x/y)$. In this case sum over $m$ first, then, similarly to above,
\begin{equation}
    \begin{aligned}
    S(k_1,k_2) & \le \sum_{\substack{y<p \le z \\ \Omega \left(p-1,t\right)=k_1}}  f(p-1) \sum_{\substack{m\le x/p \\ \Omega\left(m,t\right)=k_2}} f(m)  \\
    & \ll_{b,\eta}  \frac{x}{\log (x/y)} \exp \brackets{\sum_{t<p\le x/y} \frac{f(p)}{p}} \frac{\lambda(f)^{k_2}}{k_2!} \sum_{\substack{y<p \le z \\ \Omega \left(p-1,t\right)=k_1}} \frac{f(p-1)}{p-1} \\
      &  \ll_{b,\eta}  \frac{x}{\log (x/y)} \cdot \frac{\log (z/y)}{(\log z)^2} \cdot \exp \brackets{\sum_{t<p\le x/y} \frac{f(p)}{p} + \sum_{t<p\le z} \frac{f(p)}{p}}  \frac{\lambda(f)^{k_1}\lambda(f)^{k_2}}{k_1!k_2!}\,.
    \end{aligned}
\end{equation}
Observe that in the first case there is a missing factor of size $\log (z/y)/\log z^\prime$, however, it happens to be of order $1$ since $\log (z/y) \asymp \log z^\prime$. 

Combining the two cases, we see that
\begin{equation}
\begin{aligned}
        S(k_1,k_2)&  \ll_{b,\eta}  \frac{xu}{\log (x/y)\log z} \cdot \exp \brackets{\sum_{t<p\le x/y} \frac{f(p)}{p} + \sum_{t<p\le z} \frac{f(p)}{p}}  \frac{\lambda(f)^{k}}{k_1!k_2!}  \\
        & \ll
    \frac{u}{\log z^\prime} \exp \brackets{\sum_{t<p\le z^\prime} \frac{f(p)}{p}}\binom{k}{k_1} \cdot \pi_k(f) \,.
\end{aligned}
\end{equation}

Recognizing this bound as proportional to the term in the binomial sum $\sum_{k_1\le k} \binom{k}{k_1} = 2^k$, we indeed verify that its main contribution comes from $k_1 \approx k_2 \approx k/2$.

Now let us bound $M_{1,k}(f)$. If $k \le k_0$, Poisson estimates hold for all $k_1$ and $k_2$, so
\begin{equation}\label{eq:M1 multiplicative 1}
    \begin{aligned}
        M_{1,k}(f) & = \sum_{k_1+k_2=k} S(k_1,k_2)  \ll_{b,\eta}  
        \frac{2^k u}{\log z^\prime} \exp \brackets{\sum_{t<p\le z^\prime}\frac{f(p)}{p}} \cdot \pi_k(f) \,.
    \end{aligned}
\end{equation}

From now on assume that $k > k_0$. Fix $a \in (1/2,1)$ such that $ (1+\varepsilon)/\rho <a < (2-\varepsilon)/\rho$. Split $M_{1,k}(f) = S_1 + S_2$ into two parts such that $S_1$ contains $k_1, k_2 < a k$, and $S_2$ is its complement. The first sum $S_1$ satisfies the same bound as \eqref{eq:M1 multiplicative 1}, so we only need to bound $S_2$.

 Observe that $k_1 \ge ak$ implies $k_2 \le (1-a)k$. Then
 \begin{equation}\label{eq:M1 f S2}
     S_2 \le \sum_{\substack{m\le x/y \\ \Omega\left(m,t\right) \le (1-a)k}}  \sum_{\substack{y<p \le x/m \\ \Omega \left(p-1,t\right)\ge ak}}  f(m(p-1)) + \sum_{\substack{m\le x/y \\ \Omega\left(m,t\right) > (1-a)k}}  \sum_{\substack{y<p \le x/m \\ \Omega \left(p-1,t\right)\le ak}}  f(m(p-1))\,.
 \end{equation}

 Let $\ell = ak = a \rho\lambda(f)$. Then by Chernoff's inequality with exponent $\eta = \ell/\lambda(f) = a\rho < 2-\varepsilon$ via \Cref{lem:Tenenbaumsum} we have
\begin{equation}
    \begin{aligned}
        \sum_{\substack{n \le N\\ \Omega(n,t) > \ell}} f(n)  \le \eta^{-\ell} \sum_{\substack{n \le N}} f(n) \eta^{\Omega(n,t)} & \ll_{b,\eta} \frac{N}{\log N} \brackets{ \sum_{t<p\le N} \frac{f(p)}{p}}\brackets{\frac{\lambda(f) e}{\ell}}^\ell \\
        & = \frac{N}{\log N} \brackets{ \sum_{t<p\le N} \frac{f(p)}{p}} e^{-\lambda(f) (Q(a\rho)-1)} 
        \,;
    \end{aligned}
\end{equation}
by \Cref{lem:Koukoulopoulossum} the sum $\sum_{{n \le N,\,\Omega(p-1,t) > \ell}} f(p-1)$ is bounded by the same quantity times an additional factor of size $1/\log N$. Applying these two founds to \eqref{eq:M1 f S2} similarly to $S_1$ we obtain 
\begin{equation}\label{eq:M1k f S2 2}
    \begin{aligned}
        S_2 
        & \ll_{b,\eta} \cdot \frac{xu} {(\log z)(\log (x/y))} \brackets{ \sum_{t<p\le x/y} \frac{f(p)}{p} + \sum_{t<p\le z} \frac{f(p)}{p}} \cdot \brackets{\frac{\lambda(f) e}{\ell}}^\ell   \brackets{\frac{\lambda(f) e}{k-\ell}}^{k-\ell}
    \end{aligned}
\end{equation}
where
\begin{equation}
    \begin{aligned}
        \brackets{\frac{\lambda(f) e}{\ell}}^\ell   \brackets{\frac{\lambda(f) e}{k-\ell}}^{k-\ell} \ll e^{-\lambda(f) (Q(\rho(1-a)) + Q(\rho a)-2)}\,.
    \end{aligned}
\end{equation}
Comparing \eqref{eq:M1k f S2 2} to \eqref{eq:M1 multiplicative 1}, and observing that $(2\lambda(f))^{k}/k! = e^{-\lambda(f) (Q(\rho) -1+\rho\log 2)} k^{-1/2}$, we see that we only need to make sure that
\begin{equation}
    \begin{aligned}
        Q((1-a)\rho) + Q(a\rho) > Q(\rho) + 1 - \rho \log 2 = 2 Q(\rho/2)\,.
    \end{aligned}
\end{equation}
This is true for any $a \in (1/2,1)$ since $Q(\rho)$ is a convex function. 
\end{proof}

In the following lemma we establish the order of magnitude for the first moment of $r(n)$ for integers $n$ with fixed $k=\Omega(n,t)$.

\begin{lemma}\label{lem: M1 fixed k}
    Let $3 \le y < z\le x$ such that $z \ge 2y$ and $z^\prime = \min\{z, x/y\} \to \infty$ as $x \to \infty$. Let $z/y \le t \le z^\prime$. Let $k$ be an integer such that $\varepsilon\log\log t \le k\leq (4-\varepsilon)\log\log t$ for some small fixed $\varepsilon>0$. Then
    \begin{equation}
        M_{1,k} = \sum_{\substack{n \le x\\\Omega(n,t) = k}} r(n)
        \asymp_\varepsilon \frac{2^k}{\log t} \cdot\frac{\log (z/y)}{\log z} \cdot \frac{x(\log\log t)^k}{k!\log t}\,.
    \end{equation}
\end{lemma}
\begin{proof}
The upper found for $M_{1,k}$ follows from \Cref{lem:first moment multiplicative} with $f\equiv 1$, so we only need to establish the matching lower bound for $k \ge \varepsilon \log\log t$. 

Let $\rho = k/\log\log t < 4-\varepsilon$. Fix $a \in (1/2,1)$ such that $a\rho < 2 -\varepsilon/2$. Let $n = (p-1)m \le x$ such that $\Omega(n,t) = \Omega(p-1,t) + \Omega(m,t) = k$. Fix $k_1 = \Omega(p-1,t)$ and $k_2 = \Omega(m,t)$, then
    \begin{equation}\label{lem: M1 fixed k-1}
        M_{1,k} \ge S_1 =:\sum_{\substack{k_1+k_2=k \\ k_i \le ak}} \sum_{\substack{y<p \le z \\ \Omega \left(p-1,t\right)=k_1}} \sum_{\substack{m\le x/p \\ \Omega\left(m,t\right)=k_2}} 1  \gg \sum_{\substack{k_1+k_2=k \\ k_i \le ak}} \sum_{\substack{y<p \le z \\ k_1 \le \Omega \left(p-1,t\right) \le k_1+2}} \sum_{\substack{m\le x/p-1 \\ k_2 \le \Omega\left(m,t\right)\le k_2+2}} 1\,.
    \end{equation}  
If $\log (x/y) \le 100 \log (z/y) \le 100\log z^\prime$, sum over $p$ first, and if $\log (x/y) \ge 100 \log (z/y)$, sum over $m$ first. Then similarly to the proof of \Cref{lem:first moment multiplicative} for $f \equiv 1$, we obtain
\begin{equation}\label{eq:M1k S1}
    S_1 
    \gg 
       \frac{x\lambda^k}{k!\log t} \cdot \frac{2^k}{\log t} \cdot \frac{\log (z/y)}{\log z} = \frac{x\lambda^k}{k!\log t} \cdot \frac{2^ku}{\log t}\,,
\end{equation}
via the lower bounds in \Cref{lem:Big-Omega-bounds} and \Cref{lem:Big-Omega-bounds-shifted} respectively.
\end{proof}

\subsection{Estimates under the barrier condition}

In this section we estimate the first and the second moments of the adjusted representation function under the barrier condition for integers with fixed $\Omega(n,t)=k$ for some $t \le z^\prime = \min\{z,x/y\}$ and $\varepsilon \log\log t \le k \le (4-\varepsilon) \log\log t$.

We change the representation function slightly by restricting it to a positive proportion of divisors. That is, for $c \in (0,1/2)$ small and fixed we define
\begin{equation}\label{def:r tilde}
    \tilde{r} (n) := \#\{(p-1)d= n: p\in (y,z],\ d\in (x/z,x/y],\ P^+(p-1)\ge z^c, P^+(d) \ge (x/y)^c\}\,.
\end{equation}

We set up our barrier as follows.

Let $J = \lceil\log\log t\rceil$ and $t_j = e^{e^{j}}$ for $1\le j \le J$. Let $\{B(j)\}_{j=1}^J$ be an arbitrary sequence of positive barrier constants. Let $C \ge 100$ be large and fixed. Then we define the barrier set
 % such that $B(j) \le k$
 \begin{equation}\label{def:barrier}
         \mathcal{B}(k) = \{n\le x:  P^+(n) \ge x^c, \ \Omega(n,t) = k,\ \Omega(n,t_j) \le B(j) +C\ \forall 1\le j \le J\}\,
 \end{equation}
depending on $x,k,t,c,C$ and $\{B(j)\}_{j=1}^J$.

We first show that the barrier condition $\{\Omega(n,t_j) \le B(j)\}$ does not affect the size of the first moment of the representation function $\tilde{r}(n)$ for ``reasonable'' constants $B(j)$. In this case, ``reasonable'' means that $\Omega(n,t_j)$ stays within the Gaussian regime around its expectation.

\begin{prop}\label{prop:firstmomentbarrier}
    Let $3 \le y \le z \le x$ such that $z \ge 2y$ and $z^\prime = \min\{z, x/y\} \to \infty$ as $x \to \infty$. Let $z/y \le t \le z^\prime$ and $J = \lceil\log\log t\rceil$. 
    
    Let $k$ be an integer such that $\varepsilon\log\log t \le k\leq (4-\varepsilon)\log\log t$ for some small fixed $\varepsilon>0$. Let $\mathcal{B}(k)$ be the barrier set defined in \eqref{def:barrier} corresponding to some $c \in (0,1/2)$, some $C\ge 100$, and barrier constants $\{B(j)\}_{j=1}^J$ such that 
    \begin{equation}
        B(j) = \frac{kj}{J} + g(\min\{j,J-j\})
    \end{equation}
    for some function $g(s) = s^{1/2} \psi(s)$, where $s^{\varepsilon} \le \psi(s) \le s^{1/2-\varepsilon} $.
% Let
% \begin{equation}
%     \tilde{r} (n) = \#\{(p-1)d= n: p\in (y,z],\ d\in (x/z,x/y],\ P^+(p-1)\ge z^c, P^+(d) \ge (x/y)^c\}\,.
% \end{equation}
Then for sufficiently large $C$ and sufficiently small $c$
    \[
    M_{1,k}(\mathcal{B})= \sum_{\substack{n\leq x}}\tilde{r}(n)\bm{1}_{n \in \mathcal{B}(k)} \ge M_{1,k}/2\,.
    \]
\end{prop}

\begin{proof}
Let us show that $M_{1,k} - M_{1,k}(\mathcal{B})\le M_{1,k}/2$. Write
\begin{equation}
\begin{aligned}
         M_{1,k} - M_{1,k}(\mathcal{B}) = \sum_{\substack{n\leq x\\
    \Omega(n,t) = k}}r(n) \bm{1}_{n \notin \mathcal{B}(k)} + \sum_{\substack{n\leq x\\
    \Omega(n,t) = k}}(r(n)-\tilde{r}(n)) \bm{1}_{n \in \mathcal{B}(k)} =: S_1 + S_2\,.
\end{aligned}
\end{equation}
Let $t_j = e^{e^j}$. We have $\bm{1}_{n \notin \mathcal{B}(k)} \le \bm{1}_{\exists j \le \lambda: \Omega(n,t_j) > B(j)} + \bm{1}_{P^+(n) \le x^c}$. By the union bound
\begin{equation}
    T_1 = \sum_{\substack{n\leq x\\
    \Omega(n,t) = k}}r(n)\bm{1}_{\exists j \le J: \Omega(n,t_j) > B(j)} \le \sum_{j\le J} \sum_{\substack{n\leq x\\
    \Omega(n,t) = k}}r(n) \bm{1}_{\Omega(n,t_j) > B(j)}\,.
\end{equation}
Apply Chernoff's inequality to each summand in $j$ with exponent 
\begin{equation}
    \eta = \eta(j) = \frac{B(j) (J-j)}{j(k-B(j))} = \frac{1 + \frac{g(j)}{j}}{1-\frac{g(j)}{J-j}} < 2\,,
\end{equation}
then by \Cref{lem:first moment multiplicative} and \Cref{lem: M1 fixed k}, there exists some $C_0 > 0$ independent of $j$ and $k$ such that
\begin{equation}
    \begin{aligned}
        T_1(j) & = \sum_{\substack{n\leq x\\
    \Omega(n,t) = k}}r(n) \bm{1}_{\Omega(n,t_j) > B(j)}\le \sum_{\substack{n\leq x\\
    \Omega(n,t) = k}}r(n) \eta^{\Omega(n,t_j)-B_j}\\
    & \le C_0 M_{1,k} \brackets{1 + \frac{(\eta-1)j}{J}}^k \eta^{-B(j)}\,.
    \end{aligned}
\end{equation}
Further, there exist some $C_1,C_2 > 0$ independent of $j$ and $g$ such that
\begin{equation}
    \begin{aligned}
        \brackets{1 + \frac{(\eta-1)j}{J}}^k \eta^{-B(j)} & = \brackets{1 + \frac{g(j)}{j}}^{-B(j)} \brackets{1 - \frac{g(j)}{J-j}}^{-(k-B(j))} \\
        & \le C_1\exp\left(-\frac{C_2  g(j)^2}{\min\{j,J-j\}}\right)\,.
    \end{aligned}
\end{equation}
Then, for $g(j) = j^{1/2}\psi(j) + C$, where $\psi(j) \ge j^\varepsilon$, 
\begin{equation}
    \begin{aligned}
        T_1 &\le C_0C_1 \sum_{j \le J} \exp\left(-\frac{C_2g(j)^2}{\min\{j,J-j\}}\right) \\
        & \le 2 C_0C_1\sum_{j \le J/2} \exp\left(-C_2( \psi(j)^2 + 2C \psi(j)/\sqrt{j} + C^2/j)\right)\,.
    \end{aligned}
\end{equation}
We can ensure that $T_1\le \frac{1}{10} M_{1,k}$ for $C$ sufficiently large in terms of $\varepsilon$, $C_0$, $C_1$ and $C_2$.

Now suppose $P^+(n) \le v = x^{c} \le x^{1/2}$. Without the loss of generality we can assume $n \ge x^{1/2}$. Use $1\le (n/\sqrt{x})^{\xi}$ for $\xi = 1/(10\log v) = 1/(10c \log x)$, then by \Cref{lem:first moment multiplicative}
\begin{equation}
    \begin{aligned}
        T_2 & = \sum_{\substack{x^{1/2} \le n\leq x\\
    \Omega(n,t) = k}}r(n) \bm{1}_{P^+(n)\le v} \le x^{-\xi/2} \sum_{\substack{x^{1/2} \le n\leq x\\
    \Omega(n,t) = k}}r(n) n^{\xi} \bm{1}_{P^+(n)\le v} \\
   & \le C_3 \frac{2^ku}{\log t} \frac{x^{1-\xi/2}}{k!\log x} \brackets{\sum_{p\le v}\frac{1}{p^{1-\xi}}}^k \exp\brackets{\sum_{t<p\le v}\frac{1}{p^{1-\xi}}}\,
    \end{aligned}
\end{equation}
for some constant $C_3 > 0$ independent of $v$. We have
\begin{equation}
    \begin{aligned}
        \sum_{p\le v}\frac{1}{p^{1-\xi}} - \sum_{p\le v}\frac{1}{p} \le 2\xi \sum_{p\le v} \frac{\log p}{p} \le 4\xi \log v\le 1
    \end{aligned}
\end{equation}
so
\begin{equation}
\begin{aligned}
    T_2 \le C_3 C_4 x^{-\xi}\cdot  M_{1,k}  = C_3 C_4 e^{-\frac{1}{10c}} \cdot  M_{1,k}\,.
\end{aligned}
\end{equation}
To bound the contribution of $P^+(p-1)\le z^c$ and $P^+(d)\le (x/y)^c$ in $S_2$, one can insert the same trick into the proof of \Cref{lem:first moment multiplicative} to similarly obtain
\begin{equation}
    S_2 \le \sum_{\substack{n\leq x\\
    \Omega(n,t) = k}}(r(n)-\tilde{r}(n))  \le C_5 e^{-1/(10c)} M_{1,k} \,
\end{equation}
for some $C_5 > 0$ independent of $k,t$ and $c$.
Choosing sufficiently small $c > 0$ we can ensure that $e^{-1/(10c)}$ is as small as required in terms of $C_3$,$C_4$ and $C_5$.
\end{proof}

Next, we show that the second moment of the representation function $\tilde{r}(n)$ under the barrier condition can be bounded in terms of the barrier constants. The following proposition is not sharp when barrier constants $\{B(j)\}$ are taken to be too large, however, the proof can be easily adjusted in such a case.

\begin{prop}\label{prop:second moment general barrier}
    Let $3 \le y \le z \le x$ such that $z \ge 2y$ and $z^\prime = \min\{z, x/y\} \to \infty$ as $x \to \infty$. Let $(z^\prime)^c \le t \le z^\prime$ and $J = \lceil\log\log t\rceil$. 
    
    Let $k$ be an integer such that $\varepsilon\log\log t \le k\leq (4-\varepsilon)\log\log t$ for some small fixed $\varepsilon>0$. Let $\mathcal{B}(k)$ be the barrier set defined in \eqref{def:barrier} corresponding to some $c \in (0,1/2)$, some $C\ge 100$, and barrier constants $\{B(j)\}_{j=1}^J$ such that $B(j) \le k$.
% Let
% \begin{equation}
%     \tilde{r} (n) = \#\{(p-1)d= n: p\in (y,z],\ d\in (x/z,x/y],\ P^+(p-1)\ge z^c, P^+(d) \ge (x/y)^c\}\,.
% \end{equation}
     Then
    \[
    M_{2,k}^\prime(\mathcal{B}) = \sum_{\substack{n\leq x}}\tilde{r}(n)(\tilde{r}(n)-1)\bm{1}_{n \in \mathcal{B}(k)} \ll M_{1,k} \cdot \frac{\log (z/y)}{\log z} \sum_{j\le J} 2^{B(j)} (e^{j} + \log (z/y))^{-1}\,.
    \]
\end{prop}

\begin{proof}
Let $\chi_k(n) = \bm{1}_{n \in \mathcal{B}(k)}$. Let $P^+(n)\ge x^c$ be the largest prime factor of $n$. 

Let $n=(p-1)d$, then either $P^+(n)\mid p-1$, or $P^+(n)\mid d$, so that
\begin{equation}
\begin{aligned}
       \tilde{r}(n) & 
     \le \sum_{\substack{(p-1) \mid n,\, P^+(p-1)\ge y^c\\p\in(y,z]}} \bm{1}_{P^+(n) \nmid p-1} + \sum_{\substack{(p-1)d=n, \, P^+(d) \ge (x/y)^c \\ d \in (x/z,x/y]}} \bm{1}_{P^+(n) \mid p-1} 
    =: r_1(n) + r_2(n)
\end{aligned}
\end{equation}
and $\tilde{r}(n)(\tilde{r}(n) - 1) \le 2(r_1(n)(r_1(n) - 1) + r_2(n)(r_2(n)-1))$. Thus to bound $M_{2,k}^\prime(\mathcal{B})$ it suffices to bound 
\begin{equation}
    S_1 := \sum_{n\le x} r_1(n)(r_1(n) - 1)\chi_k(n),\quad S_2 := \sum_{n\le x} r_2(n)(r_2(n) - 1) \chi_k(n)\,.
\end{equation}
The two cases are symmetric. In both cases, we will use $P^+(n)$ to move to the logarithmic density in $n$. In the first case, we are considering shifted primes on the interval $(y,z]$, while in the second, we are considering complements of shifted primes on the interval $(x/z,x/y]$.

Let us first bound $S_1$. Let $p_1 = P^+(n)$ such that $n = n^\prime p_1$, where $P^+(n^\prime) \le p_1$. In this case we have $(p-1)p_1 \le x$, so $p_1 \le x/y$. Given $p_1\ge x^c$, this is only possible when $y \le x^{1-c}$ implying in particular that $\log t \asymp \log z$. Then $r_1(n)(r_1(n) - 1) \le r_1(n^\prime)(r_1(n^\prime) -1) $ and, since $p_1 \ge x^c \ge t$ does not affect the condition $\mathcal{B}(k)$, $\chi_k(n) = \chi_k(n^\prime)$. Sum over $p_1$ and then the $z$-rough part of $n^\prime$ (by \Cref{lem:rough}) so that
\begin{equation}
    \begin{aligned}
        S_1  \le  \sum_{n^\prime \le x} r_1(n^\prime)(r_1(n^\prime)-1)\chi_k(n^\prime) \sum_{x^c \le p_1 \le\frac{x}{n\prime}} 1
        & \ll \frac{x}{\log x} \sum_{n^\prime\le x} \frac{r_1(n^\prime)(r_1(n^\prime)-1)\chi_k(n^\prime)}{n^\prime}\\
        & \ll \frac{x}{\log z} \sum_{\substack{n\le x\\P^+(n)\le z}} \frac{r_1(n)(r_1(n)-1)\chi_k(n)}{n}\,.
    \end{aligned}
\end{equation}
Let $g=(p-1,q-1)$, then $p-1 = ga$ and $q-1=gb$ for some integers $a < b$ such that $(a,b)=1$. Since $y \le ga,gb \le z$, we have $a,b \le \min\{x/y,z\}$ and $a < b < a z/y$. Without the loss of generality, assume that $P^+(b) \ge P^+(a)$. 

Let ${p_2~=~P^+((p-1)(q-1))~=~P^+(gab)}\ge z^c$. Then either $p_2 = P^+(q-1) > P^+(p-1)$, or $P^+(p-1) = P^+(q-1)$ so that $p_2 = P^+(g) \ge P^+(b) > P^+(a)$. Denote the corresponding parts of $S_1$ by $T_1$ and $T_2$ so that 
\begin{equation} \label{eq:M2 barrier S1 T1 1}
    T_1 = \frac{1}{\log z}\sum_{\substack{n\le x\\P^+(n)\le z}}\frac{R_1(n) \chi_k(n)}{n}\,, \quad T_2 = \frac{1}{\log z}\sum_{\substack{n\le x\\P^+(n)\le z}}\frac{R_2(n) \chi_k(n)}{n}\,,
\end{equation}
where $R_1(n) := \#\{p-1,q-1\mid n: y<p<q\le z,\ P^+(q-1) > P^+(p-1) \ge z^c\}$ and $R_2(n) := \#\{p-1,q-1\mid n: y<p<q\le z,\ P^+(q-1) = P^+(p-1) \ge z^c\}$. 

We will see that in the first case $r_1(n)^2$ ``decouples'' into two statistically independent copies of $r_1(n)$ producing the square of the first moment of $r_1(n)$:
\begin{equation}\label{eq:T1 final}
    \begin{aligned}
        T_1
        & \ll \brackets{\frac{2^k u}{\log t}}^2  \cdot  \frac{\brackets{\log\log t}^k }{k!\log t}\,.
    \end{aligned}
\end{equation}
As a result, treating such pairs no longer requires the barrier. 

The second case has the larger contribution of the two with
\begin{equation}\label{eq: T2 final}
\begin{aligned}
        T_2  & \ll \frac{(\log\log t + O(1))^k}{k!\log t} \cdot \frac{2^k u^2}{\log t}  \sum_{j\le J} 2^{B(j)}  (e^{j} + \log (z/y))^{-1}\,
\end{aligned} 
\end{equation}
and genuinely requires a barrier condition to control, though we will see that the barrier condition is only checked at a single moving point. Note that the bound \eqref{eq:M2 barrier S1 T1 5} for $T_1$ corresponds to the final term at $j = J$ in $T_2$, so the bound for $S_1$ follows once it is established for $T_2$.

Let us first bound $T_1$.

Suppose $p_2 = P^+(q-1) > P^+(p-1)$, then $p_2 \le \min\{z,x/y\}$. Write $n = mp_2$ and $q-1 = dp_2$. Then $\chi_k(n) \le \bm{1}_{\Omega(m,t)\in \{k,k-1\}}$ and
\begin{equation} \label{eq:M2 barrier S1 T1 2}
    \begin{aligned}
       R_1(n) =  R_1(mp_2) \le  \sum_{\substack{p-1\mid m\\ y<p \le z}} 1 \cdot \sum_{\substack{d|m \\ y <dp_2\le z}} \bm{1}_{dp_2 + 1 \in \P} \le r(m) \cdot \sum_{\substack{d|m\\ y <dp_2\le z}} \bm{1}_{dp_2 + 1 \in \P}\,.
    \end{aligned}
\end{equation}

Split $p_2$ into $e$-adic intervals such that $p_2 \in I_i = (e^{i-1},e^{i}]$ for $c\log z \le i \le \log t$, then by \Cref{lem:Fordsieve} 
\begin{equation} \label{eq:M2 barrier S1 T1 4}
    \begin{aligned}
    \sum_{p_2 \ge z^c} \frac{1}{p_2} \sum_{\substack{d|m\\ \frac{y}{p_2} <d\le \frac{z}{p_2}}} \bm{1}_{dp_2 + 1 \in \P} & \ll \sum_{i\ge c\log z}   \sum_{\substack{d \mid m\\ \frac{y}{e^i} <  d \le \frac{4z}{e^i}}}\sum_{\substack{p_2 \in I_i\\ dp_2 + 1 \in \P}}\frac{1}{p_2} \ll \sum_{i\ge c\log z}  \frac{1}{i^2} \sum_{\substack{d \mid m\\ \frac{y}{e^i} <  d \le \frac{4z}{e^i}}} \frac{d}{\varphi(d)}\\
        & \ll \frac{1}{(\log z)^2}\sum_{\substack{d \mid m}} \frac{d}{\varphi(d)} \sum_{\log y/d \le i \le \log 4z/d} 1\\
        & \ll \frac{\log (z/y)}{(\log z)^2}\sum_{\substack{d \mid m}} \frac{d}{\varphi(d)} \ll \frac{\log (z/y)}{(\log z)^2} \cdot  \frac{\tau(m)m}{\varphi(m)} \,.
        % & \ll  e^{-2j} \min\{e^j,\log (z/y)\} f_1(m_j)  \,,
    \end{aligned}
\end{equation}
Plugging this into $T_1$, we upper-bound $T_1$ by the first (logarithmic) moment of the representation function $r(m)$ twisted by a multiplicative function:
\begin{equation}\label{eq:M2 barrier S1 T1 5}
    \begin{aligned}
        T_1 & \ll \frac{1}{\log z}\sum_{\substack{m\le x,\, P^+(m)\le z\\\Omega(m, t) =k}}  \sum_{p_2 \ge z^c} \frac{R_1(mp_2)}{mp_2} \ll 
        & \frac{\log (z/y)}{(\log z)^3} \sum_{\substack{m\le x,\, P^+(m)\le z\\\Omega(m, t) =k}} \frac{r(m)}{m} \cdot \frac{\tau(m)m}{\varphi(m)} \,.
    \end{aligned}
\end{equation}
On average $n/\varphi(n) = \prod_{p\mid n} (1-1/p)^{-1}\asymp 1$ and $\tau(m) \asymp 2^k$, though this is of course not true pointwise.  In any case, both functions trivially belong to $\mathcal{F}_{b,1}$ which allows us to apply \Cref{lem:first moment multiplicative} so that
\begin{equation}
    \begin{aligned}
        \sum_{\substack{m\le x,\, P^+(m)\le z\\\Omega(m, t) =k}} \frac{r(m)}{m} \cdot \frac{\tau(m)m}{\varphi(m)} & \ll \frac{2^k u}{\log t}  \cdot  \frac{1}{k!} \brackets{2\sum_{p\le t} \frac{1}{p-1}}^k \exp\brackets{2 \sum_{t < p \le z} \frac{1}{p-1}}\\
        & \ll \frac{2^k u}{\log t}  \cdot  \frac{2^k(\log z)^2}{\log t} \cdot \frac{\brackets{\log\log t + O(1)}^k }{k!\log t}\,.
    \end{aligned}
\end{equation}
Plugging the above into \eqref{eq:M2 barrier S1 T1 5} yields \eqref{eq:T1 final} as claimed.

We now move on to $T_2$.

Suppose $p_2 = P^+(p-1)=P^+(q-1)$. Write $n= mp_2$ and $g = g^\prime p_2$ as well as $p-1 =  g^\prime ap_2$ and $q-1= g^\prime bp_2$. Then
\begin{equation}
    \begin{aligned}
       R_2(n) =  R_2(mp_2) \le \sum_{\substack{g^\prime a b\mid m,\, (a,b)=1\\ y<g^\prime a p_2 \le z\\ a < b < az/y}}  \bm{1}_{g^\prime a p_2 + 1,g^\prime b p_2 +1 \in \P}\,.
    \end{aligned}
\end{equation}
Fix $p_2 \in I_i \in (e^{i-1},e^{i}]$, then by \Cref{lem:Fordsieve} for each $c\log z \le i \le \log t$
\begin{equation}\label{eq:M2 barrier S1 T2 0}
    \sum_{\substack{p_2 \in I_i}} \frac{1}{p_2} \cdot  \bm{1}_{g^\prime a p_2 + 1,g^\prime b p_2+1 \in \P} \ll \frac{ \Phi(g^\prime a,g^\prime b)}{i^3}\ll \frac{ \Phi(g^\prime a,g^\prime b)}{(\log z)^3}\,,
\end{equation}
where
\begin{equation}
\begin{aligned}
        \Phi(g^\prime a,g^\prime b) & = \frac{(g^\prime)^3ab(b-a)}{\varphi((g^\prime)^3ab(b-a))}\cdot\frac{g^\prime}{\varphi(g^\prime)} = \frac{g^\prime ab(b-a)}{\varphi(g^\prime ab(b-a))}\cdot\frac{g^\prime}{\varphi(g^\prime)} \,. 
\end{aligned}
\end{equation}
Observe that $g^\prime a b \mid m$, then
\begin{equation}
    \Phi(g^\prime a,g^\prime b) \le \frac{g^\prime ab}{\varphi(g^\prime ab)}\cdot \frac{b-a}{\varphi(b-a)} \cdot \frac{g^\prime}{\varphi(g^\prime)} \le \brackets{\frac{m}{\varphi(m)}}^2  \frac{b-a}{\varphi(b-a)}\,,
\end{equation}
which is much easier to manage.

Summing over $i \ge c\log z$ and noting that $\log y/(g^\prime a) \le i\le \log 4z/(g^\prime a)$ for each fixed $g^\prime, a$ we obtain
\begin{equation}\label{eq:M2 barrier S1 T2 1}
    \begin{aligned}
   \sum_{i \ge c\log z} \sum_{p_2 \in I_i} \frac{R_2(mp_2)}{p_2} & \ll \frac{\log (z/y)}{(\log z)^4} \cdot\brackets{\frac{m}{\varphi(m)}}^2\sum_{\substack{g^\prime a b\mid m,\, (a,b)=1\\ g^\prime \le z,\, a,b \le t\\ a < b < az/y }}
\frac{b-a}{\varphi(b-a)} \,.
\end{aligned}
\end{equation}
Observe that the sum over $g^\prime$ can be bounded by
\begin{equation}
    \sum_{\substack{g^\prime a b\mid m,\\ g^\prime \le z }} 1 \le \sum_{\substack{g^\prime \mid \frac{m}{ab}\\ g^\prime \le z}} 1 \le \tau(m/ab)\,,
\end{equation}
so we can rewrite the final sum in \eqref{eq:M2 barrier S1 T2 1} as
\begin{equation}
    \begin{aligned}
        R_3(m) := \sum_{\substack{g^\prime a b\mid m,\, (a,b)=1\\ g^\prime \le z,\, a,b \le t\\ a < b < az/y }}
\frac{b-a}{\varphi(b-a)} 
& \le \sum_{\substack{a b\mid m,\, (a,b)=1\\ a,b \le t\\ a < b < az/y }} \frac{b-a}{\varphi(b-a)}\cdot \tau(m/ab)\,,
    \end{aligned}
\end{equation}
which now looks quite similar to the sum over generic divisors of $m$ except for the factor
\begin{equation}\label{eq:(b-a)/phi(b-a)}
    \frac{b-a}{\varphi(b-a)}
\end{equation}
posing a minor inconvenience. We remove this factor for $(a,b)=1$ by fixing $p_3 = P^+(b)$ and bounding the number of primes $p_3$ in arithmetic progressions. 

Let $p_3 = P^+(b) \in I_i$ for some $i\le \log t$, and $i \in I_j$ for some $j\le J = \lceil \log\log t\rceil$. That is, setting $t_j = e^{e^{j}}$, we put $p_3$ in the interval $(t_{j-1},t_j] = (e^{e^{j-1}},e^{e^j}]$. Denote the $t_j$-smooth part of $m^\prime$ by $m_j^\prime = \prod_{p^\ell \| m^\prime,\, p \le t_j} p^\ell$. 

Write $\frac{n}{\varphi(n)} = \sum_{d\mid n} \frac{\mu^2(d)}{\varphi(d)}$, then
\begin{equation}
    \frac{b-a}{\varphi(b-a)} = \frac{b^\prime p_3-a}{\varphi(b^\prime p_3-a)} = \sum_{d \mid b^\prime p_3-a} \frac{\mu^2(d)}{\varphi(d)} = \sum_{d \le 2t} \frac{\mu^2(d)}{\varphi(d)} \cdot \bm{1}_{b^\prime p_3 \equiv a (d)}\,.
\end{equation}
We have $(a,b^\prime)=1$, so for each $i \in I_j$ by the Brun--Titchmarsh inequality
\begin{equation} \label{eq:M2 barrier S1 T2 2}
    \begin{aligned}
        \sum_{p_3\in I_i}\frac{1}{p_3} \cdot \frac{b^\prime p_3-a}{\varphi(b^\prime p_3-a)}&  \le \sum_{d\le 2t} \frac{\mu(d)^2}{\varphi(d)} \sum_{\substack{p_3\in I_i\\ b^\prime p_3 \equiv a (d)}} \frac{1}{p_3} \ll \sum_{d\le 2t} \frac{\mu^2(d)}{\varphi^2(d)} \cdot \frac{1}{\log (e^i/d + 2)}\\
        & \ll \sum_{d\le 2t} \frac{1}{d^{2-\varepsilon}} \cdot \frac{1}{\log (e^i/d + 2)} \ll \frac{1}{i} \ll e^{-j}\,.
    \end{aligned}
\end{equation}
The contribution from $p_3 \in (t_{j-1},t_j]$ is then bounded by
\begin{equation}
\begin{aligned}
        \sum_{p_3\in (t_{j-1},t_j]} \frac{R_3(m^\prime p_3)}{p_3} & \ll \sum_{i\in I_j} \sum_{\substack{a b^\prime \mid m^\prime_j,\, (a,b^\prime)=1\\ \frac{a}{e^i} < b^\prime < \frac{4az}{ye^i} }}  \tau(m^\prime/ab^\prime) \sum_{p_3\in I_i} \frac{1}{p_3}  \cdot \frac{b^\prime p_3-a}{\varphi(b^\prime p_3-a)}\\
         & \ll   e^{-j} \sum_{i\in I_j} \sum_{\substack{a b^\prime \mid m^\prime_j\\ \frac{a}{e^i} < b^\prime < \frac{4az}{ye^i} }}  \tau(m^\prime/ab^\prime) \\ 
          & \ll   e^{-j} \min\{e^j,\log (z/y)\} \sum_{\substack{a b^\prime \mid m^\prime_j }}  \tau(m^\prime/ab^\prime)\,.
\end{aligned}
\end{equation}
Observe that
\begin{equation} \label{eq:M2 barrier S1 T2 4}
    \begin{aligned}
        \sum_{\substack{a b^\prime \mid m^\prime_j }}  \tau(m^\prime/ab^\prime) =  \sum_{\substack{d \mid m^\prime_j}}  \tau(m^\prime/d) \tau(d) = \tau_4(m^\prime_j) \tau(m^\prime/m^\prime_j) = \frac{\tau_4(m^\prime_j)}{\tau(m^\prime_j)} \cdot \tau(m^\prime)\,.
    \end{aligned}
\end{equation}

Thus, summing over $j\le J$ and sieving the $t$-rough part (via \Cref{lem:Tenenbaumsum}) we reduce $T_2$ to
\begin{equation}
    \begin{aligned}
        T_2 & \ll  \frac{u^2}{(\log z)^2} \sum_{j\le J}\frac{1}{e^{j} + \log (z/y)} \sum_{\substack{m^\prime\le x\\ P^+(m^\prime)\le z}} \frac{\chi_k(m^\prime)\tau_4(m^\prime_j) \tau(m^\prime/m^\prime_j)}{m^\prime} \brackets{\frac{m^\prime}{\varphi(m^\prime)}}^2\\
        & \ll  \frac{u^2}{(\log t)^2} \sum_{j\le J}  \frac{1}{e^{j} + \log (z/y)}\sum_{\substack{m \le x\\ P^+(m )\le t}} \frac{\chi_k(m)\tau_4(m_j) \tau(m/m_j)}{m} \brackets{\frac{m}{\varphi(m)}}^2
        \,,
    \end{aligned}
\end{equation}
where for fixed $j$ it suffices to check the barrier condition only at $t_j$. Given $\Omega(m, t_j) \le B(j)$, we essentially have
\begin{equation}
    \tau_4(m_j) \cdot \tau(m/m_j) \le 4^{\Omega(m,t_j)} \cdot 2^{\Omega(m) - \Omega(m,t_j)} = 2^{k + \Omega(m,t_j)} \le 2^k \cdot 2^{B(j)}
\end{equation}
and thus 
\begin{equation}
\begin{aligned}
        T_2 & \ll \frac{2^k u^2}{(\log t)^2} \sum_{j\le J} 2^{B(j)}  (e^{j} + \log (z/y))^{-1} \cdot \sum_{\substack{m \le x\\ P^+(m)\le t}} \frac{\bm{1}_{\Omega(m,t) = k}}{m} \\
        & \ll \frac{2^k u^2}{\log t}  \sum_{j\le J} 2^{B(j)}  (e^{j} + \log (z/y))^{-1} \cdot  \frac{(\log\log t + O(1))^k}{k!\log t}\,.
\end{aligned} 
\end{equation}
To make the above rigorous, use Chernoff's inequality with exponent $\eta = 2$ so that $\bm{1}_{\Omega(n,t_j) \le B(j)} \le 2^{B(j) - \Omega(n,t_j)}$ and apply \Cref{lem: hardy-ramanujan multiplicative} for each $j\le J$. 

Now let us turn to the case $p_1 = P^+(n) \mid g$ and bound $S_2 = \sum_{n\le x} r_2(n)(r_2(n) -1) \chi_k(n)$. By assumption $t \le x^c \le p_1 \le z$, so this is only possible when $z\ge x^c$ and thus $\log t \asymp \log (x/y)$. If $P^+(n) \mid g$, then $P^+(n) = P^+(p-1) = P^+(q-1)$, so we will see that the contribution of this case is of the same order of magnitude as $T_2$.

Let $n = n^\prime p_1$ and $g = g^\prime p_1$, then $p-1= g^\prime p_1 a$ and $q-1 = g^\prime p_1 b$. Further, $n^\prime = g^\prime abd$ for some $d\le x/y$. Using
\begin{equation}
    r_2(n^\prime p_1) (r_2(n^\prime p_1) - 1) \le \sum_{\substack{g^\prime ab d = n^\prime,\,(a,b)=1 \\ ad<bd \in (x/z,x/y]}}\bm{1}_{g^\prime ap_1+1,g^\prime bp_1+1 \in \P}
\end{equation}
together with \Cref{lem:Fordsieve}, similarly to \eqref{eq:M2 barrier S1 T2 0} and \eqref{eq:M2 barrier S1 T2 1}, we obtain
\begin{equation}
\begin{aligned}
       \sum_{x^c \le p_1\le \frac{x}{n^\prime}} \sum_{\substack{g^\prime ab d = n^\prime,\,(a,b)=1 \\ ad<bd \in (x/z,x/y]}} \bm{1}_{g^\prime ap_1+1,g^\prime bp_1+1 \in \P} & \ll \frac{x}{n^\prime (\log x)^3} \sum_{\substack{g^\prime ab d = n^\prime,\,(a,b)=1 \\ ad<bd \in (x/z,x/y]}} \Phi(g^\prime a, g^\prime b)
\end{aligned}
\end{equation}
where
\begin{equation}
    \sum_{\substack{g^\prime ab d = n^\prime,\,(a,b)=1 \\ ad<bd \in (x/z,x/y]}} \Phi(g^\prime a, g^\prime b) \ll \brackets{\frac{n^\prime}{\varphi(n^\prime)}}^2 \sum_{\substack{g^\prime ab d = n^\prime,\,(a,b)=1 \\ ad<bd \in (x/z,x/y]}} \frac{b-a}{\varphi(b-a)}\,.
\end{equation}
Thus, we have reduced $S_2$ to
\begin{equation}
    \begin{aligned}
      S_2 & \ll \frac{x}{(\log x)^3} \sum_{n^\prime \le x} \frac{\chi_k(n^\prime)}{n^\prime} \brackets{\frac{n^\prime}{\varphi(n^\prime)}}^2 \sum_{\substack{ab d \mid  n^\prime,\,(a,b)=1 \\ ad<bd \in (x/z,x/y]}} \frac{b-a}{\varphi(b-a)}  \,.
    \end{aligned}
\end{equation}
The rest of the proof for $S_2$ is identical to one for $T_2$ with $(y,z]$ replaced by $(x/z,x/y]$. 
% The statement of the Proposition then follows.
\end{proof}

\begin{remark}
The choice of the exponent $\eta = 2$ is not necessarily optimal and can be adjusted depending on $B(j)$ to $\eta(j) = \max \{1, 2/ h(j)\}$, where
\[
h(j) =  \frac{B(j)(J-j)}{j(k-B(j))}\,.
\]    
Then $\sum_{j\le J}2^{B(j)} e^{-j}$ is replaced by
\begin{equation}
    \sum_{j\le J} \min\{2^{B(j)}, (1 + j/J)^k\} e^{-j}\,.
\end{equation}
\end{remark}